\documentclass[11pt,a4paper]{amsart}
\usepackage{amssymb, amscd, amsmath, enumerate,verbatim}

\numberwithin{equation}{subsection}
\newtheorem{theorem}{Theorem}[subsection]
\newtheorem{lemma}[theorem]{Lemma}
\newtheorem{proposition}[theorem]{Proposition}
\newtheorem{corollary}[theorem]{Corollary}

\theoremstyle{definition}

\newtheorem{definition}[theorem]{Definition}
\newtheorem{remark}[theorem]{Remark}
\newtheorem{example}[theorem]{Example}

\newcommand{\Ql}{\overline{\mathbb Q_{\l}}}
\newcommand{\mk}{\mathbf{k}}

\newcommand{\C}{\mathcal{C}}
\newcommand{\CC}{\mathbf{C}}
\newcommand{\M}{\mathcal{M}}
\newcommand{\CKM}{\mathbf{K\ddot{a}h}}
\newcommand{\Dif}{\mathbf{Dif}}
\newcommand{\Top}{\mathbf{Top}}
\newcommand{\Var}{\mathbf{V}(\mathbb{C})}
\newcommand{\Vman}{\mathbf{DMV}(\mathbb{C})}
\newcommand{\Ab}{\mathcal{A}}
\newcommand{\Chains}[1]{\mathbf{C}_*(#1)}
\newcommand{\Chainsab}{\Chains{\Ab}}
\newcommand{\sMod}{\mathbf{\Sigma Mod}}
\newcommand{\sModplus}{\mathbf{\Sigma^+ Mod}}
\newcommand{\MMod}{\mathbf{MMod}}
\newcommand{\MOp}{\mathbf{MOp}}
\newcommand{\Op}{\mathbf{Op}}
\newcommand{\Opplus}{\mathbf{Op^+}}
\newcommand{\Opn}[1]{\Op (\leq\! #1)}

\newcommand{\Aut}{\mathrm{Aut}}

\newcommand{\Hom}{\mathrm{Hom}}

\newcommand{\GL}[2]{\mathbf{GL}_{#1}({\mathbf{#2}})}
\newcommand{\AUT}[1]{\mathbf{Aut} (#1)\,}

\newcommand{\eg}[1]{\mathbf{F}_{#1}}
\newcommand{\egprima}[1]{\mathbf{F}'_{#1}}

\newcommand{\clausura}{\overline{\mathbf{k}}}

\newcommand{\invlim}{\underset{\leftarrow}{\lim}\,}
\newcommand{\alt}{\mathrm{Alt}\,}

\newcommand{\modulizero}{\overline{\mathcal{M}}_0}

\newcommand {\lra}{\longrightarrow}
\newcommand {\lla}{\longleftarrow}

\begin{document}
\begin{large}

\title[Moduli spaces and formal operads]{Moduli spaces and formal operads}

\author[F. Guill\'{e}n]{F. Guill{\'e}n Santos}
\address[F. Guill{\'e}n Santos and V. Navarro]{ Departament
d'\`{A}lgebra i Geometria\\  Universitat de Barcelona\\ Gran Via 585, 08007 Barcelona
(Spain)}
\author[V. Navarro] {V. Navarro }

\author[P. Pascual]{P. Pascual}
\address[P. Pascual and A. Roig]{ Departament de Matem\`{a}tica Aplicada
I\\ Universitat Polit\`{e}cnica de Catalunya\\Diagonal 647, 08028
Barcelona (Spain). }
\author[A. Roig]{A. Roig}
\email{guillen@mat.ub.es\\pere.pascual@upc.es\\agustin.roig@upc.es }

\maketitle

\begin{abstract}
Let  $\overline{\M}_{g,l}$ be the moduli space  of stable
algebraic curves of genus $g$ with $l$ marked points. With the
operations which relate the different moduli spaces
 identifying marked points, the family $(\overline{\M}_{g,l})_{g,l}$
is a modular operad of projective smooth
   Deligne-Mumford
stacks, $\overline{\M}$. In this paper we prove that the modular operad of singular
chains $C_*(\overline{\M}_{};\mathbb{Q})$ is formal; so it is weakly equivalent to the
modular operad of its homology $H_*(\overline{\M}_{};\mathbb{Q})$. As a consequence, the
\lq\lq up to homotopy" algebras of these two operads are the same. To obtain this result
we prove a formality theorem for operads analogous to Deligne-Griffiths-Morgan-Sullivan
formality theorem, the  existence of  minimal models of modular operads, and a
characterization of formality for operads which shows that formality is independent of
the ground field.

\end{abstract}

\section{Introduction}

In recent years, moduli spaces of Riemann surfaces such as the moduli spaces of stable
algebraic curves of genus $g$ with $l$ marked points, $\overline{\M}_{g,l}$, have played
an important role in the mathematical formulation of certain theories inspired by
physics, such as the complete cohomological field theories.

In these developments, the  operations which relate the
different moduli spaces $\overline{\M}_{g,l}$ identifying
marked points, $ \overline{\M}_{g,l} \times
\overline{\M}_{h,m} \longrightarrow \overline{\M}_{g+h,l+m-2} $ and $ \overline{\M}_{g,l}
\longrightarrow \overline{\M}_{g+1,l-2} $, have been interpreted  in terms  of operads.
With these operations  the spaces $\overline{\M}_{0,l}$,
$l\geq 3$, form a cyclic operad of projective smooth
varieties, $\overline{\M}_0$ (\cite{GeK94}), and the spaces
$\overline{\M}_{g,l}$, $g,l\geq 0$, $2g-2+l > 0$, form a
modular operad of projective smooth Deligne-Mumford stacks,
$\overline{\M}$ (\cite{GeK98}). Therefore, the homologies of
these operads, $H_*(\overline{\M}_{0};\mathbb{Q})$ and
$H_*(\overline{\M};\mathbb{Q})$, are cyclic and modular
operads respectively.

An important result in the algebraic theory of the Gromov-Witten invariants is that, if
$X$ is a complex projective manifold and $\Lambda(X)$ is the Novikov ring of $X$, the
cohomology $H^*(X;\Lambda(X))$ has a natural structure of an algebra over the modular
operad $H_*(\overline{\M}_{};\mathbb{Q})$, and so it is a complete cohomological field
theory (\cite{Be}, see \cite{Man}).

But there is another modular operad associated to the geometric operad
$\overline{\M}_{}$: the modular operad $C_*(\overline{\M}_{};\mathbb{Q})$ of singular
chains. Algebras over this operad have been studied in \cite{GeK98}, \cite{KSV} and
\cite{KVZ}.

In this paper we prove that the modular operad
$C_*(\overline{\M}_{};\mathbb{Q})$ is formal; so it is weakly
equivalent to the modular operad of its homology
$H_*(\overline{\M}_{};\mathbb{Q})$. As a consequence, the
\lq\lq up to homotopy" algebras of these two operads are the
same.

A paradigmatic example of  operad is the little $2$-disc operad of Boardman-Vogt,
$\mathcal{D}_2(l)$, of configurations of $l$ disjoint discs in the unity disc of
$\mathbb{R}^2$. Our result can be seen as the analogue for $\overline{\M}_{}$ of the
Kontsevich-Tamarkin's formality theorem of  for $C_*(\mathcal{D}_2;\mathbb{Q})$
(\cite{Ko} and \cite{T}; moreover \cite{Ko}  also explains the relation between this
formality theorem, Deligne's conjecture in Hochschild cohomology and  Kontsevich's
formality theorem in deformation quantization).

Our paper is organized as  follows. In section 2, we study
symmetric monoidal functors between symmetric monoidal
categories, since they induce functors between the categories
of their operads. After recalling some definitions and fixing
some notations of operads and monoidal categories, we prove a
symmetric De Rham theorem.  We then introduce the notion of
formal symmetric monoidal functor, and we see how this kind of
functor produces formal operads.

In section 3, as a consequence of Hodge theory, we prove that
the cubic chain functor on the category of compact K\"{a}hler
manifolds $C_*: \CKM \longrightarrow \Chains{\mathbb{R}}$ is a
formal symmetric monoidal  functor. It follows that, if $X$ is
an operad of compact K\"{a}hler manifolds, then the operad of
chains $C_*(X;\Bbb R)$ is formal. This is the analogue in the
theory of operads of the Deligne-Griffiths-Morgan-Sullivan
formality theorem in rational homotopy theory (\cite{DGMS}).

The goal of  sections  $4$, $5$ and $6$ is to prove the { descent\/} of formality from
$\mathbb{R}$ to $\mathbb{Q}$. In section $4$ we recall some results due to M. Markl on
minimal models of operads in the form that we will use in order to generalize them to
cyclic and modular operads.

In section 5, drawing on Deligne's weight theory for Frobenius
endomorphism in \'{e}tale cohomology,  we introduce weights and
show the formality of the category of complexes endowed with a
pure endomorphism. Next, in  th. \ref{aixeca}, we prove a
characterization of formality of an operad in terms of the
lifting of automorphisms of the homology of the operad to
automorphisms of the operad itself.

The automorphism group of a minimal operad with homology of
finite type is a pro-algebraic group. This result allows us to
use the descent theory of algebraic groups to prove the
independence of formality of the ground field in th.
\ref{descens}.

In section 7 we show how  the above  results can be extended
easily to cyclic operads. In particular we obtain the
formality of the cyclic operad  $C_*(\modulizero ; \mathbb{Q}
)$.

In the last section, we go one step further and  prove  the
above  results also for modular operads. In particular, we
introduce minimal models of modular operads and we prove their
existence and lifting properties. Here, we follow
Grothendieck's idea in his ``jeu de L\'{e}go-Teichm\"{u}ller"
(\cite{Gro}), in which he builds the complete Teichm\"{u}ller
tower inductively on the modular dimension.
 Once this is  established, the proofs of
the previous sections can be transferred to the modular context without difficulty.
Finally, we conclude the formality of the modular operad $C_*(\overline{\mathcal M}_{};
\mathbb{Q})$.

\section{Formal operads}

\subsection{Operads} Let us recall some definitions and notations about operads (see
\cite{GK}, \cite{KM}, \cite{MSS}).

\subsubsection{}Let $\Sigma$ be the {\it symmetric groupoid\/}, that is,
the category whose objects are the  sets $\{1,\dots , n \} $,
$n\ge 1$ , and the only morphisms are those of  the symmetric
groups $\Sigma_n = Aut\{1,\dots , n \} $.

\subsubsection{} Let $\mathcal C$ be a  category. The category of
contravariant functors from $\Sigma$ to $\mathcal{C}$ is
called the category of $\Sigma$-{\it modules\/} and is denoted
by $\sMod_{\mathcal{C}}$, or just $\sMod$ if $\mathcal C$ is
understood. We identify its objects with sequences of objects
in $\mathcal{C}$, $E=\left((E(l)\right)_{l\geq 1}$, with a
right $\Sigma_l$-action on each $E(l)$. If $e$ is an element
of $E(l)$, $l$ is called the {\it arity\/} of $e$. If $E$ and
$F$ are $\Sigma$-modules, a {\it morphism of $\Sigma$-modules}
$f:E \longrightarrow F$ is a sequence of $\Sigma_l$-{
equivariant\/} morphisms $f(l): E(l) \longrightarrow F(l), \
l\geq 1$.

\subsubsection{}Let $(\mathcal{C}, \otimes, \mathbf{1})$ be a symmetric monoidal
category. A {\it unital $\Sigma$-operad\/} (an {\it operad\/}
for short) in $\mathcal{C}$ is a $\Sigma$-module $ P$ together
with a family of { structure morphisms\/}: { composition}
$\gamma_{l;m_1,\dots ,m_l} : P(l)\otimes P(m_1) \otimes \dots
\otimes P(m_l)
\longrightarrow P(m_1 + \cdots + m_l)$, and { unit\/} $\eta : \mathbf{1} \longrightarrow
P(1)$, satisfying the axioms of equivariance, associativity,
and unit. A {\it morphism of operads\/} is a morphism of
$\Sigma$-modules compatible with structure morphisms. Let us
denote by $\Op_{\mathcal{C}}$, or simply $\Op$ when $\mathcal
C$ is understood, the category of operads in $\mathcal C$ and
its morphisms.

\subsection{Symmetric monoidal categories and  functors}
In the study of $\Sigma$-operads the commutativity constraint
plays an important role. In particular the  functors we are
interested in are functors between symmetric monoidal
categories which are  compatible with the associativity,
commutativity and unit constraints.

\subsubsection{}
The following are some of the  symmetric monoidal categories
we will deal with in this paper. On the one hand, the
geometric ones:

$\Top $: the category of topological spaces.

$\Dif$: the category of differentiable manifolds.

$\CKM$: the category of compact K\"{a}hler manifolds.

$\Var $: the category of smooth projective $\mathbb
C$-schemes.

On the other hand, the algebraic categories, which will be
subcategories, or variants of

$\Chainsab $: the category of   complexes with a differential
of degree $-1$ of an abelian monoidal symmetric category $(\Ab
, \otimes , \mathbf{1})$. The morphisms are called chain maps.
If $\Ab$ is the category of $R$-modules for some ring $R$, we
will denote it by $\Chains{R}$. Operads in $\Chainsab$ are
also called {\it dg operads\/}.

In a symmetric monoidal category $(\mathcal C,\otimes ,\mathbf
1)$ we usually denote the natural commutativity isomorphism by
$\tau_{X,Y}:X\otimes Y\longrightarrow Y\otimes X$. For
example, in $\Chainsab $,  the natural commutativity
isomorphism
$$
\tau_{X,Y} : X \otimes Y \longrightarrow Y \otimes X
$$
includes the signs:
$$
\tau_{X,Y} (x \otimes y) = (-1)^{\deg (x)\deg (y)} y \otimes x\ .$$

\subsubsection{}
As usual, we move from a geometric category to an algebraic
one through a  functor.
 Let us recall (see \cite{KS}) that a {\it  monoidal functor\/}
$$
(F , \kappa ,\eta ) : (\mathcal{C}, \otimes , \mathbf{1})
\longrightarrow (\mathcal{D},
\otimes , \mathbf{1}')
$$
between monoidal  categories is a functor $F : \mathcal{C}
\longrightarrow \mathcal{D}$ together with a natural morphism
of $\mathcal{D}$, $$ \kappa_{X,Y} : FX \otimes FY
\longrightarrow F(X \otimes Y),$$ for all objects $X, Y \in \mathcal{C}$, and a morphism
of $\mathcal{D}$, $ \eta : \mathbf{1}' \longrightarrow
F\mathbf{1}$, compatibles with the constraints of
associativity, and unit. We will refer the {\it K\"{u}nneth
morphism} as $\kappa$.

 If $\mathcal C$ and $\mathcal D$ are
symmetric monoidal categories, a monoidal functor $F :
\mathcal{C} \longrightarrow
\mathcal{D}$ is said to be {\it symmetric} if $\kappa$ is compatible with the
commutativity constraint.

For example, the {\it homology functor\/} $H_*: \Chainsab
\longrightarrow \Chainsab$ is a symmetric monoidal functor,
taking the usual K{\"u}nneth morphism
$$H_*(X)\otimes H_*(Y) \longrightarrow H_*(X \otimes Y).$$
as $\kappa$.

 Let $F,G:\mathcal{C}\rightrightarrows
\mathcal{D}$ be two monoidal functors. A natural transformation $ \phi : F \Rightarrow G
$ is said to be {\it monoidal\/} if it is compatible with
$\kappa$ and $\eta$.

\subsubsection{}
Let $F : \mathcal{C} \longrightarrow \mathcal{D}$ be a
symmetric monoidal functor. It is easy to prove that, applied
componentwise, $F$ induces a functor between $\Sigma$-operads
$$\Op_{F}:\Op_{\mathcal{C}} \longrightarrow \Op_{\mathcal{D}} \ , $$ also denoted by $F$.

 In
particular,  for an operad $P\in \Op_{\Chainsab}$, its
homology is an operad $HP \in
\Op_{\Chainsab}$.

In the same way, if $F, G:\mathcal{C} \rightrightarrows
\mathcal{D}$ are two symmetric monoidal functors, a monoidal
natural transformation $\phi:F\Rightarrow G$ induces a natural
transformation
$$\Op_\phi:\Op_{F}\Rightarrow \Op_G,$$ also
denoted  by $\phi$.

\subsection{Weak equivalences} We will use  weak equivalences in
several contexts.

\subsubsection{} Let $X$ and $Y$ be objects of $\CC_*(\Ab)$. A chain map $f:
X \longrightarrow Y$ is said to be a {\it weak equivalence of
 complexes \/} if the induced morphism $f_* = Hf : HX
\longrightarrow HY$ is an isomorphism.

\subsubsection{} Let $\mathcal C$ be a category,  $\Ab$ an
abelian category, and  $F,G:\mathcal C\rightrightarrows
\mathbf C_*(\Ab)$ two functors. A natural transformation
$\phi:F\Rightarrow G $ is said to be a {\it weak equivalence
of functors}, if the morphism $\phi(X):F(X)\rightarrow G(X)$
is a weak equivalence, for every object $X$ in $\mathcal C$.

\subsubsection{} A
morphism $\rho : P \longrightarrow Q $ of operads in
$\Chainsab$ is said to be a {\it weak equivalence
 of operads \/} if  $\rho (l) : P(l)
\longrightarrow Q(l)$ is a  weak equivalence of chain complexes, for all $l$.

\subsubsection{} Let  $\mathcal C$ be a  category endowed with a distinguished
 class of morphism
called {\it weak equivalences}. We suppose that this is a
saturated class of morphisms which contains  all isomorphisms.
Two objects $X$ and $Y$ of $\mathcal C$ are said to be {\it
weakly equivalent} if there
 exists a sequence of morphism of $\mathcal C$
$$
X \longleftarrow X_1 \longrightarrow \cdots \longleftarrow
X_{n-1} \longrightarrow Y.
$$
which are weak equivalences. If $X$ and $Y$ are weakly
equivalent, we say that $Y$ is a {\it model} of $X$.

\subsubsection{}
The following proposition  is an easy consequence of the definitions.
\begin{proposition}\label{Operadsequivalents}
If $F,G:\mathcal C\rightrightarrows \mathcal D$ are two weakly
equivalent  symmetric monoidal functors, the functors $\Op_F$
and $\Op_G$ are weakly equivalent. \! In particular, for every
operad $P$ in $\mathcal C$, the operads $F(P)$ and $G(P)$ are
weakly equivalent.
\end{proposition}

\subsection{Symmetric De Rham's theorem}
In this section we will demonstrate a symmetric version of  De
Rham's theorem comparing the complex of singular cubic chains
and  De Rham's complex of currents, including its symmetric
structure as monoidal functors.

\subsubsection{} It is well known that
the functor  of {\it    singular chains \/}$$ S_*(\phantom{X}
; \mathbb{Z}) : \Top
\longrightarrow \Chains{\mathbb{Z}} \ ,
$$
together with the {\it shuffle product\/}
 as the K\"{u}nneth morphism, is a symmetric monoidal functor.

On the other hand, let $\mathcal D'_*(M)$ be the complex of De
Rham's {\it currents\/} of a differentiable manifold $M$; that
is, $\mathcal D'_*(M)$ is the topological dual of the complex
$\mathcal D^*(M)$ of differential forms with compact support.
Then the functor
$$
\mathcal{D}'_* : \Dif \longrightarrow \Chains{\mathbb{R}}  \ ,
$$
is a symmetric monoidal functor  with the K\"{u}nneth morphism
$$
\kappa_{M,N}: \mathcal{D}'_*(M) \otimes \mathcal{D}'_*(N) \longrightarrow
\mathcal{D}'_*(M \times N)
$$
induced by the tensor product of currents. Thereby, if $S \in
\mathcal{D}'_*(M)$, and $T
\in \mathcal{D}'_*(N)$, then
$$
<\!\kappa(S \otimes T)\,,  \,\pi_M^*(\omega ) \wedge \pi_N^*
(\nu )\!> \,=\, <\! S\,,\,\omega \!>\,\cdot \,<\!T\,,\,\nu\! >
$$
for all $\omega \in \mathcal{D}^*(M)$ and $\nu \in
\mathcal{D}^*(N)$.

 In order to compare the functor of
currents with the functor of  singular chains on
differentiable manifolds, one can consider the complex of
chains $S_*^\infty(M)$ generated by the $\mathcal
C^{\infty}$-maps $\Delta^p \longrightarrow M$. The
corresponding functor of {\it $\mathcal C^{\infty}$-singular
chains\/} $S_*^\infty:\Dif \longrightarrow\bold C _*(\mathbb
Z)$ is also a symmetric monoidal functor with the shuffle
product.

On the one hand, the natural inclusion of $\mathcal
C^{\infty}$-singular chains
 into
singular ones defines a monoidal natural transformation
$S_*^\infty \Rightarrow S_* :\Dif
\rightrightarrows \bf C_*(\mathbb Z), $ and from the approximation theorem it follows
that it   is a
  weak equivalence of functors.
On the other hand, by Stokes' theorem, integration along
$\C^\infty$-singular simplexes induces a natural
transformation $ \int:S_*^\infty
\Rightarrow\mathcal{D}'_*:\Dif
\rightrightarrows \bold C_*(\mathbb R) $, which is a weak equivalence of functors by De
Rham's theorem.

 So the functors of singular  chains
and currents are weakly equivalent. However, the natural transformation $\int$  is not
compatible with the monoidal structures. To overcome this deficiency, we will use the
cubic  singular chains  instead of the simplicial ones.

\subsubsection{} For a topological space $X$, cubic
chains are generated by the singular cubes of $X$; that is,
the continuous maps $I^p\longrightarrow X$, where $I$ is the
unit interval of the real line $\mathbb R$ and $p\in \mathbb
N$, modulo the degenerate ones  (see e.g. \cite{Mas}). The
functor of {\it cubic\/} chains
$$
C_*(\phantom{X} ; \mathbb{Z}) : \Top \longrightarrow
\Chains{\mathbb{Z}} \ ,
$$
 with the  {\it cross product\/}
$$
\times\;: C_*(X ; \mathbb{Z}) \otimes C_*(Y ; \mathbb{Z})
 \longrightarrow C_*(X \times Y ; \mathbb{Z}) \ ,
$$
which for singular cubes $c : I^p \longrightarrow X$ and $d :
I^q \longrightarrow Y$ is defined as the cartesian product
$$
c \times d : I^{p+q} = I^p \times I^q \longrightarrow X \times
Y\, ,
$$
is a monoidal functor.

If $M$ is a differentiable manifold, one may consider chains
generated by $\C^\infty$-maps $ I^p \longrightarrow M$. The
corresponding functor of {\it $\C^\infty$-cubic chains\/}
$C_*^\infty(\phantom{X} ; \mathbb{Z}):\Dif
\longrightarrow\bold C _*(\mathbb Z)$ is also a monoidal functor.

Finally the inclusion $C^{\infty}_*\Rightarrow C_*:\Dif\longrightarrow \CC_*(\mathbb Z)$
and integration $\int:C_*^{\infty}\Rightarrow \mathcal D'_*:\Dif \longrightarrow \mathbf
C_*(\mathbb R)$
 are  monoidal natural transformations,
 which are weak equivalences of monoidal functors.

 In spite of this, in this case another problem arises:
the monoidal functor of cubic chains  is not symmetric,
because the cross product is not commutative.

\subsubsection{}
However, over $\mathbb{Q}$, or more generally over a $\mathbb
Q$-algebra, we can symmetrize cross product with the classical
alternating operator. For every $n$, define the map $ \alt :
C_n(X;\mathbb{Q}) \longrightarrow C_n(X;\mathbb{Q}) $ as
$$
\alt (c) = \frac{1}{n!} \sum_{\sigma \in \Sigma_n} (-1)^{|\sigma|} c \circ \sigma \ ,
$$
for all singular cubes $c:I^p \longrightarrow X$. The map
$\alt$
 is well defined because, if $c$ is a degenerate singular cube, then $\alt(c)$
is a degenerate cubic chain.

\begin{proposition}
The functor of cubic chains $C_*(\phantom{X} ;
\mathbb{Q}):\Top \longrightarrow
\Chains{\mathbb{Q}} $, together with the { product\/}
$$
 \kappa_{X,Y} : C_*(X;\mathbb{Q}) \otimes C_*(Y;\mathbb{Q})
  \longrightarrow C_*(X \times Y;\mathbb{Q})
 \,,
 $$
defined by $\kappa(c\otimes d)=\alt(c\times d)$, is a
symmetric monoidal functor.
\end{proposition}

\begin{proof}
First of all, compatibility with units is trivial. Next,
 compatibility
 with   associativity and  commutativity constraints  is proved
  in a similar way to the usual proofs of the
associativity and commutativity properties for the wedge
product of multi-linear alternating tensors. Finally, the fact
that $\alt$ is a chain map is possibly a classical result, but
we have not found a proof in the literature. So  let us see
that $\alt$ is a chain map. Recall the definition of the
  differential
 $d:C_n(X;\mathbb{Q}) \longrightarrow C_{n-1}(X;\mathbb{Q})$ of the cubic chain
 complex.
For  $1\le i\le n,\ \epsilon\in\{0,1\}$, let
$\delta_i^{\epsilon}:I^{n-1} \longrightarrow I^{n}$ denote the
face defined by $ \delta^{\epsilon}_{i}
(t_1,\dots,t_{n-1})=(t_1,\dots,\epsilon,\dots,t_{n-1})\,\,, $
 where $\epsilon$ is in
the $i$-th place. Now, if $c\in C_n(X,\mathbb Q)$, $d(c)$
 is defined by
$$
d (c)=\sum_{i,\epsilon}(-1)^{i+\epsilon}c\circ
\delta_i^{\epsilon}\,\, ,
$$
and we have
$$
d(\alt (c)) =\frac{1}{n!}\sum_{i,\epsilon}\sum_{\sigma\in
\Sigma_n}(-1)^{i+\epsilon+ |\sigma|} c \circ \sigma\circ
\delta_i^{\epsilon} \ ,$$ and $$ \alt(d c)   =
\frac{1}{(n-1)!}\sum_{\tau\in \Sigma_{n-1}}\sum_{r,\epsilon}(-1)^{|\tau |+r+\epsilon}c
\circ \delta_r^{\epsilon}\circ\tau
$$
Then, it is easy to check the following claim

{\bf Claim:} {\it For all $\tau\in \Sigma_{n-1},$ and
$r,i\in\{1,\dots,n\}$ let
$$
\sigma_{\tau,r,i}:=(r,\dots,n)\circ \overline{\tau}\circ (i,\dots,n)^{-1}\,,
$$
where $\overline{\tau}\in \Sigma_{n}$ is defined by $
\overline{\tau}(i)=  \tau(i)$ for all $ 1\leq i < n $, and
  $(r,\dots,n)$ denotes the cycle
  $$(1,2,\dots, r-1,r,r+1,\dots,n)\mapsto (1,2\dots,r-1,r+1,\dots,n,r)\ .$$
 Then $\sigma_{\tau,r,i}$ is the only
permutation that satisfies $ \sigma_{\tau,r,i}\circ
\delta^\epsilon_i=\delta^\epsilon_r\circ \tau. $

Moreover, the set map  $$\begin{aligned} \Sigma_{n-1}\times
\{1,\dots,n\}\times
\{1,\dots,n\} &\longrightarrow \Sigma_n\times \{1,\dots,n\}\\(\tau,r,i)&\mapsto
(\sigma_{\tau,r,i},i)\end{aligned}$$ is bijective. }

To finish checking that $ d\circ \alt= \alt\circ d \  $, it
suffices to  note that, from $
(-1)^{\vert\sigma_{\tau,r,i}\vert+i}=(-1)^{|\tau|+r} \ , $ it
follows
$$
\frac{1}{n}\sum_{i=1}^{n}(-1)^{ |\sigma |+i+\epsilon}c \circ \sigma_{\tau,r,i}\circ
\delta_i^{\epsilon}=(-1)^{|\tau | +r+\epsilon}c \circ \delta_r^{\epsilon}\circ\tau
$$
for all $\epsilon$, $r$, $\tau$.
\end{proof}

The  functor of {\it $\mathcal C^{\infty}$-cubic  chains\/}
$C_*^\infty(\phantom{X};\mathbb Q):\Dif \longrightarrow\bold C
_*(\mathbb Q)$ is also a symmetric monoidal functor with the
corresponding product.

We will now   check that integration is compatible with the
monoidal structure.
\begin{lemma} Integration along $C^\infty$-cubes induces
a monoidal natural transformation
$$
\int:C_*^\infty \Rightarrow\mathcal{D}'_*:\Dif \rightrightarrows \bold C_*(\mathbb R) ,
$$
which is a weak equivalence of symmetric monoidal functors.
\end{lemma}
\begin{proof} Let $M$ be a differentiable manifold and $c: I^p \longrightarrow
M$ a $\mathcal C^\infty$-singular cube. Integration along $c$,
 $\int^M_c
\omega=\int_{I^{p}}c^*(\omega)$, defines a chain map
$$
\int^M :C_*^\infty(M; \mathbb{R}) \longrightarrow \mathcal{D}'_*(M;\mathbb{R})
$$
 by $c \mapsto  \int_c^M  $, which is obviously natural
in $M$, and a weak equivalence  by  De Rham's theorem.

The natural transformation $\int^M$ is monoidal, by Fubini's
theorem  and by the change of variables theorems. Indeed, for
$c:I^p \longrightarrow M$ and $d: I^q \longrightarrow N$
$\mathcal C^\infty$-singular cubes, and   $\omega \in
\mathcal{D}^p(M)$ and $\nu \in
\mathcal{D}^q(N)$, we have
\begin{align} \int^{M\times
N}_{\kappa  (c\otimes d)} \pi^*_M (\omega) \wedge
\pi^*_N(\nu)\notag
              &=\frac{1}{(p+q)!} \sum_{\sigma
              \in
             \Sigma_{p+q}}(-1)^{|\sigma |}
             \int_{(c\times d)\circ \sigma}^{M\times N} \pi^*_M (\omega)
             \wedge \pi^*_N(\nu)        \notag\\
              &=  \int_{c\times d}^{M\times N}\! \pi^*_M (\omega)
              \wedge \pi^*_N(\nu)     \notag\\
             & =   \left(\int_c^M \pi^*_M (\omega ) \right)
             \left(\int_{d}^N \pi^*_N (\nu ) \right)      \notag\\
             & =  \left(\kappa \circ \left(\int^M _c\otimes \int^N _d\right)
             \right)(\pi^*_M (\omega) \wedge
             \pi^*_N(\nu))\ .\notag
\end{align}
\end{proof}
To sum up, we  can state the following symmetric version of
De Rham's theorem
\begin{theorem}\label{cubicchainswecurrents}
The functors   of cubic  chains and currents
$$C_*,\mathcal
D'_*:\Dif\rightrightarrows \CC_*(\mathbb R)
 $$
are weakly equivalent symmetric monoidal functors.
\end{theorem}
\begin{remark}
A similar result can be obtained with  the functor  of {\it
oriented cubic chains \/} $C_*^{\mathrm{or}}$, used by
Kontsevich (see \cite{Ko}, (2.2)),  which, with the cross
product, is  a monoidal symmetric functor from $\Top$ to
$\Chains{\mathbb{Z}}$. It can be proved using the alternating
operator that, at least over $\mathbb{Q}$,  computes the usual
homology of a topological space. This solution is equivalent
to the previous one, because the natural projection $
C_*(\phantom{X};\mathbb Q) \longrightarrow
C_*^{\mathrm{or}}(\phantom{X};\mathbb Q)$ is a monoidal
natural transformation, which is a weak equivalence of
symmetric monoidal functors.
\end{remark}

\begin{remark} The
 above results can be seen as  analogous in the
covariant setting to the following well known statements. The
contravariant  functor of singular cochains
$S^*:\Top^{op}\longrightarrow \CC^*(\mathbb Z)$, together with
the Alexander-Whitney product, is a monoidal contravariant
functor. Although the functor of $\C^\infty$-singular cochains
$S^*_\infty:\Dif^{op}\lra \CC^*(\mathbb R)$ is weakly
equivalent to the  symmetric monoidal functor of  differential
forms $\mathcal E^*:\Dif^{op}\longrightarrow \CC^*(\mathbb
R)$, the monoidal functor $S^*_\infty$  is symmetric only up
to homotopy.  $\mathcal E^*$ can be topologically defined as a
symmetric monoidal functor using the Sullivan's $\mathbb
Q$-cdga $Su_\mathbb Q$ of simplicial differential forms, which
defines a symmetric monoidal contravariant functor
$\Top^{op}\longrightarrow \CC^*(\mathbb Q)$, together with  a
weak equivalence of symmetric monoidal  functors $\mathcal
E^*\Rightarrow Su_\mathbb R:\Dif\longrightarrow
\CC^*(\mathbb R)$.

\end{remark}

\subsection{Formality}

\subsubsection{}
The notion of formality has attracted  interest since Sullivan's work  on rational
homotopy theory. In the operadic setting the notion of formality appears in \cite{Mkl}
and \cite{Ko}.

\begin{definition} An operad $P$ in $\Chainsab$ is said to be {\it formal\/}
if it is weakly equivalent to its homology $HP$.
\end{definition}

More generally, we can give the following definition
\begin{definition}
Let $\mathcal C$ be a category endowed with an idempotent
endofunctor $H:\mathcal C\lra
\mathcal C$,   and take as  weak equivalences the morphisms  $f:X\lra Y$ such that $H(f)$
is an isomorphism. An object $X$ of $\mathcal C$ is said to be
{\it formal} if $X$ and $HX$ are weakly equivalent.
\end{definition}

\subsubsection{} In particular, a  functor $F: \mathcal{C} \longrightarrow
\Chainsab$ is  {\it formal} if it is weakly equivalent to its homology $HF: \mathcal{C}
\longrightarrow \Chainsab$. However, we will use this notion in the context of symmetric
monoidal functors. So,  the definition of formality in this
case is

\begin{definition} Let $\mathcal C$ be a symmetric monoidal category,
 and $F: \mathcal{C} \longrightarrow
\Chainsab$ a  symmetric monoidal functor. It is said that $F$ is a  {\it formal symmetric
monoidal functor\/} if $F$ and $HF$ are weakly equivalent in
the category of symmetric monoidal functors.

If the identity functor of $\Chainsab$ is a formal symmetric
functor, $\Chainsab$ is said to be  a {\it formal symmetric
monoidal category}.

Let $\mathcal B$ be a symmetric monoidal subcategory of
$\Chainsab$. If the inclusion functor  $\mathcal B\lra
\Chainsab$  is a formal symmetric monoidal functor, we will
say that $\mathcal B$ is a {\it formal symmetric monoidal
subcategory of $\Chainsab$}.

\end{definition}

The  properties below follow immediately from the definitions.

\begin{proposition}
Let $\mathcal B$ be  a  formal symmetric monoidal subcategory
of $\Chainsab$. If $F:\C\longrightarrow \Chainsab$ is a
monoidal functor with  values in the subcategory $\mathcal B$,
then $F$  is a formal symmetric monoidal functor.
\end{proposition}
\begin{proposition}\label{formalfunctorscreateformaloperads} Let $F: \mathcal{C}
\longrightarrow \Chainsab$ be a functor. If $F$ is a formal symmetric monoidal functor,
then
$$
F: \Op_{\mathcal{C}} \longrightarrow \Op_{\Chainsab}
$$
sends operads in $\mathcal{C}$ to formal operads in
$\Chainsab$.
\end{proposition}

\subsubsection{} Let $R$ be a commutative ring, and
$R-\mathbf{cdga}$  the category of differential
graded-commutative $R$-algebras, or simply cdg $R$-algebras.
It is a symmetric monoidal category. Then, if $\mathcal C$ a
symmetric monoidal category and $F: \mathcal{C}
\longrightarrow \mathbf C^*(R)$ is a symmetric monoidal
functor, it is a well known fact that $F$ induces a functor
from the category of commutative monoids of $(\mathcal
C,\otimes,\mathbf 1)$ to cdg $R$-algebras.

Besides, if $F$ is a formal symmetric monoidal functor,
 then $F$ sends commutative monoids
to formal cdg $R$-algebras.

 If $\mathcal C$ is a category with  finite
products, and  a final object $\mathbf 1$, then $(\mathcal C,\times,\mathbf 1)$, is a
symmetric monoidal category. In this case every object $X$ of $\mathcal C$ is a comonoid
object with the diagonal $X\rightarrow X\times X$ and the unit $X\rightarrow \mathbf 1$.
So we have
\begin{proposition}\label{formalidaddealgebras}
Let  $\mathcal C$ be a category with  finite products, and  a
final object $\mathbf 1$. Every  formal symmetric monoidal
contravariant functor $F:\mathcal{C}^{op}
\longrightarrow \mathbf C^*(R) $ sends objects in $\mathcal C$ to formal cdg
$R$-algebras.
\end{proposition}

\section{Hodge theory implies formality}

In  \cite{DGMS}, Deligne et al. prove the formality of the De Rham cdg algebra of a
compact K{\"a}hler manifold. In this section we will see how this result can be mimicked
for the cubic chain complex of an operad of  compact K{\"a}hler manifolds.
\subsection{Formality of De Rham's functor}
 In {\cite{DGMS}}, the first of the proofs of formality
 th. 5.22
relies on the {\it Hodge decomposition\/} for the complex of
forms and the K{\"a}hler identities. From them, the
$dd^c$-{\it lemma\/}  is proved and the existence of a diagram
of complexes, called {\it $d^c$-diagram},
$$
(\mathcal{E}^*(M),d) \longleftarrow \ (^c\mathcal{E}^*(M),d)
\longrightarrow (H^*_{d^c}(M),d) \ ,
$$
 is deduced. Here $M$ is a compact K{\"a}hler manifold, $\mathcal{E}^*(M)$ is the real De
Rham complex of $M$, $\ ^c\mathcal{E}^*(M)$ the subcomplex of
$d^c$-closed forms, and $H_{d^c}^*(M)$ the quotient complex $\
^c\mathcal{E}^*(M) / d^c(\mathcal{E}^*(M))$. In the
$d^c$-diagram, both maps are weak equivalences of chain
complexes, and the differential induced by $d$ on
$H_{d^c}^*(M)$ is zero. Since $\ ^c\mathcal{E}^*$ is a
symmetric monoidal subfunctor of $\mathcal{E}^*$ and the
morphisms of this diagram are natural, the functor
$\mathcal{E}^*$ is formal.
 So the theorem of formality can also be
stated with the previous definitions as follows

\begin{theorem} The functor of differential forms
$\mathcal{E}^* : \CKM^{op} \longrightarrow \mathbf
C^*(\mathbb{R})$ is a formal symmetric monoidal
 functor.
\end{theorem}

This result, together with prop. \ref{formalidaddealgebras}, implies the formality
theorem for  De Rham's cdg-algebra in its usual formulation: The De Rham functor
$\mathcal{E}^*: \CKM^{op} \longrightarrow \mathbb R -\mathbf{cdga}$ sends objects in
$\CKM$ to formal cdg $\mathbb R$-algebras.

\subsection{Formality of the  current complex functor}
We claim that an analogous theorem of formality  is obtained
replacing forms with currents.

\begin{theorem}\label{curretsformalfunctor} The functor of currents
$\mathcal{D}'_* : \CKM \longrightarrow \Chains{\mathbb{R}}$ is
a formal symmetric monoidal functor.
\end{theorem}

\begin{proof}
Let $M$ be a compact K\"{a}hler manifold.  It is a classical result of Hodge theory (see
\cite{Sch}) that the K{\"a}hler identities between the operators $d, d^c, \Delta ,\dots $
of the De Rham complex of differential forms are also satisfied by the corresponding dual
operators on  De Rham  complex of currents. Hence we have the following $dd^c$-lemma.

\begin{lemma} Let $T$ be a $d^c$-closed and $d$-exact  current.
Then, there exists a current $S$ such that $T = dd^cS$.
\end{lemma}

From this lemma, we can follow verbatim the first proof of
theorem 5.22 ({\it the $d^c$-Diagram Method\/}) in \cite{DGMS}
and we obtain a $d^c$-diagram for currents:
$$
(\mathcal{D}'_*(M),d) \longleftarrow \ (^c\mathcal{D}'_*(M),d)
\longrightarrow (H^{d^c}_*(M),d)  \ .
$$
Here $\ ^c\mathcal{D}'_*(M)$ denotes the subcomplex of
$\mathcal D'_*(M)$ defined by the $d^c$-closed currents, and
$H^{d^c}_*(M)$ is the quotient $\ ^c\mathcal{D}'_*(M) / d^c
(\mathcal{D}'_*(M))$. In this $d^c$-diagram both maps are weak
equivalences,  and the differential induced by $d$ on the
latter is zero. So we have $H_*^{d^c}(M)\cong H_*(\mathcal
D'_*(M)$.

Now, since $d^c$ satisfies the Leibnitz rule, $\
^c\mathcal{D}'_*$ is a symmetric monoidal subfunctor of
$\mathcal D'_*$ .

Finally, since  the morphisms of the above $d^c$-diagram are
natural and compatible with the K\"{u}nneth morphism,  it follows
that $\mathcal{D}'_* $ is a formal symmetric monoidal functor.
\end{proof}
As a consequence of the formality of the current functor and the symmetric De Rham
theorem for currents (th. \ref{cubicchainswecurrents}),  the formality of the cubic
chains functor for compact K\"{a}hler manifolds follows.
\begin{corollary}\label{cadenesformals} The functor   of cubic chains
$C_*(\phantom{X}; \mathbb R): \CKM \longrightarrow
\Chains{\mathbb{R}}$ is a formal symmetric monoidal  functor.
\end{corollary}

\subsection{Formality of K\"{a}hlerian operads}\label{hodge}

From \ref{formalfunctorscreateformaloperads} and
\ref{cadenesformals} we obtain the operadic version  of the
formality DGMS theorem.
\begin{theorem}\label{cadenesformalsR}
If $X$ is an operad in $\CKM$, then the operad of cubic chains
$C_*({X}; \mathbb R)$ is formal.
\end{theorem}

This result, together with  prop.
\ref{formalfunctorscreateformaloperads} and
 the descent theorem  th. \ref{descens}
below, implies the formality of the operad of cubic chains
with rational coefficients for every operad of compact K\"{a}hler
manifolds (see cor. \ref{cadenesformalsQ} below).

\subsection{Formality of $DM$-operads}
 The above results can be easily generalized
to the category of {\it Deligne-Mumford projective and smooth
stacks\/} over $\mathbb C$, which we will denote by $\Vman$.

Indeed, every  stack of this kind defines a compact K\"{a}hler $V$-manifold and for such
$V$-manifolds we have the functors of cubic chains, $\C^\infty$-cubic chains and
currents, and also Hodge theory  (see \cite{Ba}). This  allows us to obtain an analogous
result to cor. \ref{cadenesformals}:

\begin{theorem}\label{hodgeforVman} The functor of cubic chains
 $C_*(\phantom{X};\mathbb R): \Vman \longrightarrow \Chains{\mathbb{R}}$
is a formal symmetric monoidal  functor.
\end{theorem}

And, from \ref{formalfunctorscreateformaloperads}  follows
\begin{theorem}\label{Coperadformal} If $X$ is an
operad in $\Vman$, then the operad of cubic chains
$C_*(X;\mathbb{R})$ is formal.
\end{theorem}

\section{Minimal operads}

In this section $\mathbf{k}$ will denote a field of
characteristic zero, and  an operad will be an operad in the
category of dg vector spaces over $\mathbf {k}$,
$\CC_*(\mathbf k)$. The category of these operads is denoted
simply by $\Op$. It  is a complete and cocomplete category
(see \cite{Hi}).

\subsection{Some preliminaries}
Let us start by recalling some basic results on minimal
operads due to M. Markl (\cite{Mkl}, see \cite{MSS}).

\subsubsection{}
A {\it minimal operad\/} is an  operad of the form $(\Gamma
(V),d_M)$, where $\Gamma:\sMod\lra\Op$ is the free operad
functor, $V$ is a $\Sigma$-module with zero differential with
$V(1)=0$, and the differential $d_M$ is {\it decomposable\/}.

The free operad functor $\Gamma:\sMod\lra\Op$ is a right
adjoint functor for the forgetful functor $U:\Op\lra \sMod$.

A {\it minimal model\/} of an operad $P$ is a minimal operad
$P_\infty$, together with a weak equivalence $P_{\infty}
\longrightarrow P$.


Let $P=\left(P(l)\right)_{l\ge 1}$ be an operad.  M. Markl has
proved that, if $HP(1) =
\mathbf{k}$,  $P$  has a minimal model $P_\infty $ with $P_\infty(1)=\mathbf{k}$
(\cite{MSS}, th. 3.125).

As  observed in \cite{MSS}, remark II.1.62, the category of
operads $P$ with $P(1) =
\mathbf{k}$ is equivalent to the category of {\it pseudo-operads \/} $Q$ with $Q(1)=
\mathbf{0}$, the zero dg vector space (see {\it op. cit.}\, def. II.1.16).

{\it In the sequel, we will work only with pseudo-operads,
with $HP(1)=0$, and we will call them simply operads. We will
denote by $\Op$ the category of these operads, and
$$
\circ_i:P(l)\otimes P(m)\lra P(l+m-1),\quad 1\le i\le l,
$$
their composition operations}.

\subsection{Truncated operads}\label{Truncatedoperads}
We will now introduce the arity truncation and their  right and left adjoints, which
enables us to introduce in the operadic setting the analogs of the skeleton and
coskeleton functors of simplicial set theory.

 Here we establish the
results for the arity truncation in a form that can be easily
translated to modular operads in \S $8$.

\subsubsection{}\label{propiedadideal}
Let $E=\left(E(l)\right)_{l\ge 1}$ be a $\Sigma$-module, and
$n\ge 1$ an integer. The grading of $E$ induces a decreasing
filtration  $\left(E(\ge \!l)\right)_{l\ge 1}$, by the
sub-$\Sigma$-modules
$$
E(\ge\! l):=\left(E(i)\right)_{i\ge l}\, .
$$

Let $P$ be an operad.  We will denote by $P\cdot P(n)$ the
sub-$\Sigma$-module consisting of elements
$\alpha\circ_i\beta$ with $\alpha\in P(l)$ and $\beta\in
P(m)$, such that at least one of $l,m$ is $n$. It follows from
the definitions that
$$
P\cdot P(n)\subset P(\ge\! n).
$$
If, moreover, $P(1)=0$, then
$$
P\cdot P(n)\subset  P(\ge\! n\!+\!1\!).
$$
The first property implies that $P(\ge \!n)$ is an ideal of
$P$,   so the quotient $P/P(\ge\! n\!+\!1\!)$ is an operad,
which is zero in arities $>n$. This is a so-called
$n$-truncated operad. However, we find it more natural to give
the following definition of $n$-truncated operad.

\begin{definition}
A {\it $n$-truncated operad} is a finite sequence of objects
in $\CC(\mathbf k)$,
$$
P=(P(1),\dots,P(n)),
$$
with a right $\Sigma_l$-action on each $P(l)$, together with a
family of composition operations, satisfying those axioms of
composition operations in $\Op$ that make sense for truncated
operads. A {\it morphism of $n$-truncated operads} $f: P
\longrightarrow Q$ is a finite sequence of morphisms of
$\Sigma_l$-modules $f(l) : P(l) \longrightarrow Q(l) , \ 1\le
l\le  n$, which commute with composition operations.

Let $\Opn{n}$ denote the category of $n$-{\it truncated
operads\/} of $\CC_*(\mathbf k)$.
\end{definition}

A {\it weak equivalence\/} of $n$-truncated operads is a
morphism of $n$-truncated operads $\phi : P \longrightarrow Q$
which induces isomorphisms of graded $\mathbf{k}$-vector
spaces, $H\phi (l) : HP(l) \longrightarrow HQ(l)$, for $l = 1,
\dots , n$.

Given an operad $P$,  $t_nP:=(P(1),\dots, P(n))$ defines  a
{\it truncation functor\/}
$$
t_{n}:\Op \longrightarrow \Opn{n}\ .
$$

\subsubsection{}
For a $n$-truncated operad $P$ denote by  $t_*P$ the $\Sigma$-module that is $\mathbf 0$
in arities $>n$ and coincides with $P$ in arities $\le n$. Since  $P\cdot P(n)\subset
P(\ge\! n)$, $t_*P$ together with the structural morphisms of $P$ trivially extended, is
an operad, and the proposition below follows easily from the definitions
\begin{proposition} \label{propiedadest*}
Let $n\ge 1$ be an integer.  Then

\begin{enumerate}
\item $t_{*}:\Op(\le \! n)\lra \Op$ is a right adjoint functor for
$t_{n}$.

\item There exists a canonical isomorphism
$t_{n}\circ t_{*}\cong Id_{\Op(\le n)}$.

\item  $t_{*}$ is a fully faithful functor.

\item $t_{*}$ preserves limits.

\item For $m\ge n$, there exists a natural morphism
$$
\psi_{m,n}:t_*t_{m} \lra t_{*}t_{n}
$$
such that $\psi_{l,n}=\psi_{m,n}\circ\psi_{l,m}$, for $l\ge
m\ge n$. For an operad  $P$, the family $\left(t_{*}t_{
n}P\right)_{n}$, with the morphisms $\psi_{m,n}$, is an
inverse system of operads. The family of unit morphisms of the
adjunctions
$$
\psi_n:P\lra t_{*}t_{n}P,
$$
induces an isomorphism $\psi:P\lra
\lim\limits_{\leftarrow}t_{*}t_{n}P$.

\item Let $P, Q$ be operads. If $n\ge 2$, $t_{n-1}P=0$, and $Q\cong
t_{*}t_{n}Q$, then
$$\Hom_{\Op}(P,Q)\cong \Hom_{\Sigma_n}(P(n),Q(n)).$$
\end{enumerate}
\end{proposition}

\subsubsection{} On the other hand, the functor $t_{n}$
 also has a  left adjoint.
For a $n$-truncated operad $P$, denote by $t_{!}P$  the operad obtained freely adding to
$P$ the operations generated in arities $> n$, that is,
$$
t_{!}P=\Gamma(Ut_*P)/J\ ,
$$
where $J$ is the ideal in $\Gamma (Ut_*P)$ generated by the kernel of $ t_{
n}\Gamma(Ut_*P)\lra P$.
\begin{proposition}\label{propiedadest!}
Let $n\ge 1$ be an integer.  Then
\begin{enumerate}
\item $t_{!}:\Op(\le n)\lra \Op$ is a left adjoint functor for
$t_{n}$. \vskip 2mm

\item  There exists a canonical isomorphism  $t_{n}\circ
t_{!}\cong Id_{Op(\le n)}$.

\item  $t_{!}$ is a fully faithful functor.

\item $t_{!}$ preserves colimits.

\item For $m\le n$, there exists a natural morphism
$$\phi_{m,n}:t_{!}t_{m} \lra t_{!}t_{n}$$ such that
$\phi_{l,n}=\phi_{m,n}\circ\phi_{l,m}$, for $l\le m\le n$. For
an operad $P$, the family $\left(t_{!}t_{n}P\right)_{n}$, with
the morphisms $\phi_{m,n}$, is a directed system of operads.
The family of
 unit morphisms of the adjunctions
$$
\phi_n:t_{!}t_{n}P\lra P\,,
$$
induces an isomorphism $\phi:
\lim\limits_{\rightarrow}t_{!}t_{n}P\lra P$. \hskip 2mm

\item Let $P, Q$ be operads. If $t_{n-1}Q=0$,    $P\cong
t_{!}t_{n}P$, and $P(1)=0$, then
$$\Hom_{\Op}(P,Q)\cong \Hom_{\Sigma_n }(P(n),Q(n)).$$

\end{enumerate}

\end{proposition}
\begin{proof}
Part (1) follows from the definition of $t_!$, and the remaining parts follow from (1)
and prop. \ref{propiedadest*}.
\end{proof}

We will call the direct system of operads given by
$$
0\rightarrow t_{!}t_{1}P\rightarrow\cdots \rightarrow t_{!}t_{n-1} P \rightarrow
t_{!}t_{n}P\rightarrow\cdots
$$
the {\it canonical tower} of $P$.

As an easy consequence of the existence of right and left
adjoint functors for $t_{n}$ we obtain the following result.
\begin{corollary}
The truncation functors $t_{n}$ preserve limits and colimits.
In particular, they commute with homology,  send weak
equivalences to weak equivalences, and preserve formality.
\end{corollary}

\subsection{Principal extensions}
Next, we  recall the definition of a principal extension of
operads and show that the canonical tower of a minimal operad
is a sequence of principal extensions. This will allow us to
extend these notions to the truncated setting.

\subsubsection{}
To begin with, we establish  some notations on suspension and mapping cones of complexes
in an additive category.

If $A$ is a chain complex and $n$ is an integer, we denote by
$A[n]$ the complex defined by $A[n]_i = A_{i-n}$ with the
differential given by $d_{A[n]}=(-1)^n d_A$.

For a chain map $\eta: B\longrightarrow A$ we will denote by
$C\eta$, or by $A\oplus_\eta B[1]$, the mapping cone of
$\eta$, that is to say, the complex that in degree $i$ is
given by $(C\eta)_i=A_i\oplus B_{i-1}$ with the differential
$d(a,b)=(d_Aa +\eta b, -d_Bb)$. Therefore $ C\eta$ comes with
a canonical chain map $i_A:A\lra C\eta$ and a canonical
homogeneous map of graded objects $j_B:B[1]\lra C\eta$.

For a chain complex $X$, a chain map $\phi :
C\eta\longrightarrow X$ is determined by the chain map $\phi
i_A:A\longrightarrow X$ together with the homogeneous map
$\phi j_B :B[1]\longrightarrow X$. Conversely, if $f:A\lra X$
is a chain map and $g:B[1]\lra X$ is a homogeneous map such
that $f\eta = d_Xg + gd_B$, that is, $g$ is a homotopy between
$f\eta$ and $0$, then there exists a unique chain map
$\phi:C\eta\lra X$ such that $\phi i_A=f$ and $\phi j_B=g$. In
other words, $C\eta$ represents the functor
 $h_{\eta}:\CC_*(\mathbf k)\lra \mathbf{Sets}$ defined, for $X\in
\CC_*(\mathbf k)$, by
$$
h_{\eta}(X)=  \{(f,g); \ f\in \Hom_{\CC_*(\mathbf k)}(A, X)\ ,\ \, g\in \Hom_{\mathbf
k}(B, X)_{1}\ , \ d_Xg+gd_B=f\eta\ \},
$$
where $\Hom_{\mathbf k}(B, X)_1$  denotes the set of homogeneous maps of degree $1$ of
graded $\mathbf k$-vectorial spaces.

\subsubsection{} Recall the construction of standard cofibrations
introduced  in (\cite{Hi}). Let $P$ be an operad, $V$ a dg
$\Sigma$-module and $\xi :V[-1] \longrightarrow P$ a chain map
of dg $\Sigma$-modules. The standard cofibration associated to
these data, denoted by $P\!\langle\!V,\xi\!\rangle $ in
\cite{Hi}, is an operad that represents the functor
$h_{\xi}:\Op\lra \mathbf{Sets}$ defined, for $Q\in
\Op$, by
$$
h_{\xi}(Q)=  \{(f,g);\ f\in \Hom_{\Op}(P, Q)\ ,\, g\in \Hom_{\bf Gr\Sigma Mod}(V,
UQ)_{0}\ , \ d_Qg-gd_V=f\xi\ \},
$$
where $\Hom_{\bf Gr\Sigma Mod}(V, UQ))_0$  denotes the set of homogeneous maps of degree
$0$ of graded $\Sigma$-modules. When $V$ has zero differential, this construction is
called a {\it principal extension} and denoted by $P*_\xi \Gamma(V)$, see \cite{MSS}. For
reasons that will become clear at once, we will denote it by $ P\sqcup_\xi V$. From the
definition it follows that $ P\sqcup_\xi V$ comes with a canonical morphism of operads
$i_P:P\lra P\sqcup_\xi V$ and a canonical homogeneous map of degree $0$ of graded
$\Sigma$-modules $j_V:V\lra P\sqcup_\xi V$.

Now, one can express $P\sqcup_\xi V$ as a push-out. Let $C(V[-1])$ be the mapping cone of
id$_{V[-1]}$, $S(V)=\Gamma(V[-1])$, $T(V)=\Gamma(C(V[-1]))$ , $i_V:S(V)\lra T(V)$ the
morphism of operads induced by the canonical chain map $i:V[-1]\lra C(V[-1])$, and
$\widetilde\xi:S(V)\lra P$ the morphism of operads induced by $\xi$. Then $P\sqcup_\xi V$
is isomorphic to the push-out of the following diagram of operads
$$
  \begin{CD}
    S(V) @>\widetilde\xi>> P \\
    @Vi_{V}VV@.\\
    T(V)\ .@.
  \end{CD}
$$
If $V$ is concentrated in arity $n$, and its differential is
zero, the operad $P\sqcup_\xi V$ is called  an {\it arity $n$
principal extension}.

Let us explicitly describe it in the case that $n\ge 2$,  $P(1)=0$, and $P\cong
t_{!}t_{n}P$. First of all, since, for a truncated operad $Q\in \Op(\le n-1)$, there is a
chain of isomorphisms
$$
\begin{aligned}\label{}
\Hom(t_{n-1}(P\sqcup_\xi V),Q)&\cong \Hom(P\sqcup_\xi V, t_{*}Q)
\\
  &\cong \Hom(P,t_{*}Q)
\\
  &\cong \Hom(t_{n-1}P,Q),
\end{aligned}
$$
we have  $t_{n-1}(P\sqcup_\xi V)\cong t_{n-1}P$. Next, let $X$
be the $n$-truncated operad extending $t_{n}P$ defined by
$$
X(i)=\left\{
\aligned \label{}&P(i),\hskip 13mm \text{ if } i<n,\\
 &P( n)\oplus_\xi V,\hskip 3mm \text{ if } i=n,\endaligned
\right.
$$
the composition operations involving $V$ being trivial,
because $P(1)=0$. Then, it is clear that $X$ represents the
functor $h_{\xi}$ restricted to the category $\Op(\le n)$, so
$t_{n}\left(P\sqcup_\xi V\right)\cong X$. Finally, it is easy
to check that $t_{!}X$ satisfies the universal property of
$P\sqcup_\xi V.$ Summing up, we have proven:

\begin{proposition}
\label{estructuradeunaextensionprincipal}Let $n\ge 2$ an integer.   Let $P$ be an operad
such that $P(1)=0$ and $t_!t_{n}P\cong P$,  $V$ a dg $\Sigma$-module concentrated in
arity $n$ with zero differential, and $\xi:V[-1]\lra P(n)$ a chain map of
$\Sigma_n$-modules. The principal extension $P\sqcup_\xi V$ satisfies:

(1)\quad  $t_{n-1}(P\sqcup_\xi V)\cong t_{n-1}P$.

(2)\quad  $(P\sqcup_\xi V)(n)\cong C(\xi),$ in particular,
there  exists an exact sequence of complexes
$$
0\lra P(n)\lra (P\sqcup_\xi V)(n)\lra V\lra 0.
$$
(3) \quad $P\sqcup_\xi V\cong t_{!}t_{n}(P\sqcup_\xi V) $.

(4) A morphism of  operads $\phi:P\sqcup_\xi V\lra Q$ is
determined by a morphism of $n$-truncated  operads $f:t_nP\lra
t_nQ$, and a homogeneous map of $\Sigma_n$-modules $g:V\lra
Q(n)$, such that $f\xi =dg$.

\end{proposition}
These  results extend trivially to truncated operads.

\subsection{Minimal objects}
Now we can translate the definition of minimality of operads
of dg modules in terms of the canonical tower:

\begin{proposition}
An operad $M$ is  minimal if, and only if, $M(1)=0$ and the
canonical tower of $M$
$$
0=t_{!}t_{1}M\rightarrow\cdots \rightarrow t_{!}t_{n-1}M
\rightarrow t_{!}t_{n}M\rightarrow\cdots
$$
is a sequence of principal extensions.
\end{proposition}

\begin{proof} Let  $M = (\Gamma (V),d_M)$ be a minimal operad.
Then (see \cite{MSS}, formula II.(3.89))
$$t_{!}t_{n}M\cong(\Gamma(V(\le n)),\partial_n),$$ and
$t_{!}t_{n}M$ is an arity $n$ principal extension of
$t_{!}t_{n-1}M$ defined by $\partial_n:V(n)\lra \left(
t_{!}t_{n-1}M\right)(n)$.

Conversely, let us suppose  $M(1)=0$, and that
$t_{!}t_{n-1}M\lra t_{!}t_{n}M$ is an arity $n$ principal
extension defined by a $\Sigma$-module $V(n)$ concentrated in
arity $n$ and zero differential, for each $n$. Then,
$M=\Gamma(\bigoplus_{n\ge\! 2}V( n))$ and its  differential is
decomposable, because $M(1)=0$. So $M$ is a minimal operad.
\end{proof}

\subsubsection{}  We now give the definition of minimality for truncated
operads.

For $m\leq n$ we have an obvious truncation functor $$
t_{m}:\Opn{n} \longrightarrow
\Opn{m},$$ which has a right adjoint $t_{*}$ and a left adjoint $t_{!}$.

\begin{definition}
A $n$-truncated operad $M$ is said to be {\it minimal} if
$M(1)=0$ and the canonical tower
$$0=t_!t_{1}M\lra
t_!t_{2}M\lra \cdots \lra t_!t_{n-1}M\lra M$$ is a sequence of
($n$-truncated) principal extensions.

An operad  $M$ is said to be {\it $n$-minimal} if the
truncation $t_{n}M$ is minimal.
\end{definition}

It follows from the definitions that an operad $M$ is
$n$-minimal if, and only if, $t_{!}t_{n}M$ is minimal. It is
clear that an operad $M$ is minimal if and only if it is
$n$-minimal for every $n$, and that theorems 3.120, 3.123 and
3.125 of \cite{MSS} remain true in $\Opn{n}$, merely replacing
\lq\lq minimal"  by \lq\lq $n$-minimal".

\subsubsection{}
The category $\Op$ has a natural structure of closed model
category (\cite{Hi}). For our present purposes, we will not
need all the model structure, only a small piece: the notion
of homotopy between morphisms of operads and the fact that
minimal operads are cofibrant objects in $\Op$;  this can be
developed independently, as  in \cite{MSS}, II.3.10. From
these results  the next one follows easily.

\begin{proposition}\label{almendruco} Let $M$ be a minimal
operad and $P$ a suboperad. If the inclusion $P
\hookrightarrow M$ is a weak equivalence, then $P = M$.
\end{proposition}

\begin{proof} Let us call $i : P \hookrightarrow M$ the inclusion.
By \cite{MSS}, th. II.3.123, we can lift, up to homotopy, the
identity of $M$ in the diagram below
$$
\begin{CD}
@.                        P   \\
@.                     @VV i V \\
M    @>\mathrm{id}>>     M  \ .
\end{CD}
$$
So we obtain a morphism of operads $f : M \longrightarrow P $
such that $i f$ is homotopic to $\mathrm{id}$. Hence $if$ is a
weak equivalence and,
 by \cite{MSS}, prop.
II.3.120, it is an isomorphism. Therefore $i$ is an isomorphism too.
\end{proof}

\subsection{Automorphisms of a formal minimal operad}
For an operad $P$, let $\Aut (P)$ denote the group of its
automorphisms. The following lifting property from
automorphisms of the homology of the operad to automorphisms
of the operad itself is the first part of the characterization
of formality that we will establish in th. \ref{aixeca}.

\begin{proposition} \label{ooth} Let $M$ be a minimal operad. If $M$
is  formal, then the map $H: \Aut (M) \longrightarrow \Aut (HM)$ is surjective.
\end{proposition}

\begin{proof} Because $M$ is a formal operad, we have a sequence
of  weak equivalences
$$
M \longleftarrow X_1\lra X_2\lla\cdots \lra X_{n-1}\lla
X_{n}\lra HM \, .
$$
By the lifting property of minimal operads (\cite{MSS}, th.
II.3.123) there exists a weak equivalence
$$
\rho: M \longrightarrow HM \ .
$$
Let $\phi\in \Aut (HM)$. Again by the lifting property of
minimal operads, given the diagram

$$
\begin{CD}
@.                        M   \\
@.                     @VV \rho V \\
M    @>(H\rho) \phi (H\rho)^{-1} \rho>> HM \ .
\end{CD}
$$
there exists a morphism   $f : M \longrightarrow M$ such that
$\rho f$ is homotopic to $(H\rho) \phi (H\rho)^{-1}\rho$.
Since homotopic maps induce the same morphism in homology, it
turns out that  $f$ is a weak equivalence and, by (\cite{MSS}
II.th.3.120), it is also an isomorphism, because $M$ is
minimal. Finally, from $(H\rho)(H f )= (H\rho)\phi
(H\rho)^{-1}(H\rho)$,  $Hf= \phi$ follows.
\end{proof}

It is clear that  prop. \ref{ooth} remains true in $\Opn{n}$,
merely replacing \lq\lq minimal" and \lq\lq formal" by \lq\lq
$n$-minimal" and \lq\lq $n$-formal", respectively.

\subsection{Finiteness of the minimal model}
In this section we will show that we can transfer the
finiteness conditions of the homology of an operad to the
finiteness conditions of its minimal model.

\begin{definition}
A  $\Sigma$-module $V$ is said to be of {\it finite type} if,
for every $l$, $V(l)$  is a finite dimensional
$\mathbf{k}$-vector space. An operad $P$ is said to be {\it of
finite type\/} if the underlying $\Sigma$-module, $UP$, is of
finite type.
\end{definition}

\begin{example}
If $V$ is a $\Sigma$-module of finite type such that $V(1)=0$,
then the free operad $\Gamma(V)$ is of finite type, because
$$
\Gamma(V)(n)\cong\bigoplus_{T\in \mathcal{T}{ree}(n)} V(T),
$$
where $ \mathcal{T}{ree}(n)$ is the finite set of isomorphism
classes of $n$-labelled reduced trees, and
$V(T)=\bigotimes_{v\in\text {Vert}(T)}V(\text{In}(v))$, for
every $n$-labelled tree $T$, (see \cite{MSS}, II.1.84). In
particular, if $P$ is a $n$-truncated operad of finite type,
then $t_{!}P\cong\Gamma(P)/J$ (see prop.
\ref{propiedadest!}) is of finite type as well.
\end{example}

\begin{theorem}\label{finitudob}
Let $P$ be an operad. If the homology of $P$ is of finite
type, then every minimal model $P_{\infty}$ of $P$ is of
finite type.
\end{theorem}

\begin{proof} Let $M$ be  a minimal operad such that $HM$ is of finite type.
Since  $$M(n)=(t_{!}t_{n}M)(n),$$ it suffices to check that
$t_{!}t_{n}M$ is of finite type. We proceed by induction. The
first step of the induction is trivial because $M(1)=0$. Then,
$t_{!}t_{n}M$ is an arity $n$ principal extension of the
operad $t_{!}t_{n-1} M$ by the vector space
$$V(n)=
HC\left(\left(t_{!}t_{n-1}M\right)(n)\rightarrow
M(n)\right),$$ thus $t_{n}\left(t_{!}t_{n}M\right)$ is finite
dimensional, by the induction hypothesis. Therefore
$t_{!}t_{n}M=t_{!}\left(t_{n}t_!t_{n}M\right)$ is also of
finite type, by the previous example.
\end{proof}

\section{Weight theory implies formality}

In analogy to the formality theorem of \cite{DGMS} for the rational homotopy type of a
compact K\"{a}hler manifold, Deligne (\cite{D}) proved  formality of the ``$\Ql$-homotopy
type" of a smooth projective variety defined over a finite field using the weights of the
Frobenius action in the $\l$-adic cohomology and his solution of the Riemann hypothesis.
In this section we follow this approach and introduce weights to establish a criterion of
formality for operads based on the formality of the category of pure complexes, defined
below.

In this section $\mathbf{k}$ will denote a field of
characteristic zero, and  an operad will be an operad in
$\CC_*(\mathbf k)$.

\subsection{Weights}\label{weights}
A {\it weight function} on $\mathbf k$ is a group morphism
$w:\Gamma\longrightarrow
\mathbb Z$ defined on a subgroup $\Gamma$ of the multiplicative group  $\overline
{\mathbf k}^*$ of an algebraic closure $\overline {\mathbf k}$
of $\mathbf k$.

An element $\lambda\in \overline{ \mathbf k} $ is said to be
{\it pure of weight} $n$ if $\lambda\in \Gamma$ and
$w(\lambda)=n$. A polynomial $q(t)\in \mathbf k[t]$ is said to
be {\it pure of weight $n$} if all the roots of $q(t)$ in $
\overline{\mathbf k}$ are pure of weight $n$.
 We will write $w(q)=n$ in this case.

Let $f$ be  an endomorphism of a $\mathbf k$-vector space $V$.
If $V$ is of finite dimension, we will say that $f$ is {\it
pure of weight $n$} if its characteristic polynomial is pure
of weight $n$. When $f$ is understood,  we  will say  that $V$
is pure of weight $n$.

Let $f$ be an endomorphism of a finite dimensional $\mathbf
k$-vector space $V$. If $q(t)\in \mathbf k[t]$ is an
irreducible polynomial, we will denote by $\ker q(f)^\infty$
the primary component corresponding to the irreducible
polynomial $q(t)$, that is, the union of the subspaces $\ker
q(f)^n$, $n\ge 1$. The space $V$ decomposes as a direct sum of
primary components $V=\bigoplus \ker q(f)^\infty$, where
$q(t)$ runs through the set of all irreducible factors of the
minimal polynomial of $f$. The sum of the primary components
corresponding to the pure polynomials of weight $n$ will be
denoted by $V^n$, that is $V^n=\bigoplus_{w(q)=n}\ker
q(f)^\infty$. Hence we have a decomposition
$$
V=\bigoplus_n V^n\oplus C,
$$
where $C$ is the sum of the primary components corresponding
to the polynomials which are not pure. This decomposition will
be called the {\it weight decomposition} of $V$.

This weight decomposition is obviously functorial on the
category of pairs $(V,f)$.

Let $P$ be a complex of $\mathbf k$-vector spaces such that
$HP$ is of finite type, that is, $H_iP$ is finite dimensional,
for all $i$. An endomorphism $f$ of $P$ is said to be {\it
pure of weight $n$} if $H_i(f)$ is pure of weight $n+i$, for
all $i$. In that case, we will say that $(P,f)$, or simply
$P$, is {\it pure of weight} $n$, if the endomorphism $f$ is
understood. Obviously, if $(P,f)$ is pure of weight $n$, so do
is $(HP,Hf)$.

The following example will be useful in the sequel. Take
$\alpha\in \mathbf k^*$ that is not a root of unity and define
$w:\{\alpha^{n};\ n\in \mathbb Z\ \}\lra\mathbb Z $ by
$w(\alpha^n)=n$. Let $P$ be a finite type complex with zero
differential. Then, the {\it grading automorphism}
$\phi_\alpha$ of $P$, defined by \
$\phi_\alpha=\alpha^{i}\cdot
\mathrm{id}$\ on $P_i$, for all $i\in \mathbb Z$, is pure of weight $0$.

 Let $\mathbb F$ be a finite field of
characteristic $p$ and $q$ elements, $\l$ a prime $\not=p$. In
\cite{D}, Deligne defined a weight function $w$ on the field
$\mathbb Q_l$ as follows. If $\iota: {\Ql} \lra
\mathbb C$ is an embedding, and
$$\Gamma=\{\alpha\in \Ql\,;\, \exists n\in \mathbb Z \text { such
that } |\iota \alpha|=q^{\frac{n}{2}}\},$$ then $w:\Gamma \lra
\mathbb Z$ is defined by $w(\alpha)=n$, if $ |\iota
\alpha|=q^{\frac{n}{2}}$. The Riemann hypothesis, proved by
Deligne, asserts that the Frobenius action is a pure
endomorphism (of weight $0$) of the \'{e}tale cohomology
$H^*(X,\mathbb Q_{\l})$ of every smooth projective $\mathbb
F$-scheme $X$.

\subsection{Formality criterion}
Let $w$ be a weight function on $\mathbf k$ and denote by
$\CC^w_*(\mathbf k)$ the category of couples $(P,f)$ where $P$
is a finite type complex and $f$ is an endomorphism of $P$
which is pure of weight $0$.

\begin{theorem}\label{categoriadepesosformal}
$\CC^w_*(\mathbf k)$ is a formal symmetric monoidal category.
\end{theorem}

\begin{proof} First of all, $\mathbf 1$, with the identity, is pure of weight $0$.
Next, by the K\"{u}nneth theorem and elementary linear algebra, if
$(P,f)$ and $(Q,g)$ are pure complexes of weights $n$ and $m$
respectively, then $(P\otimes Q, f\otimes g)$ is pure of
weight $n+m$. Then it is easy to check that $\CC^w_*(\mathbf
k)$ has  a structure of symmetric monoidal category such that
the assignment
$$\CC^w_*(\mk)\lra \CC_*(\mk), \quad (P,f)\mapsto P$$ is a symmetric
monoidal functor, and that the functor of homology
$H:\CC^w_*(\mathbf k)\lra
\CC^w_*(\mathbf k)$ is a  symmetric monoidal functor as well.

To prove that $\CC^w_*(\mathbf k)$ is formal we will use the
weight function $w$ to define a symmetric monoidal functor $T$
and weak equivalences
$$
\mbox{id}_{\CC^w_*(\mathbf k)}\longleftarrow T \longrightarrow H.
$$

Let $(P,f)$ an object $\CC^w_*(\mathbf k)$. Since each $P_i$
is finite dimensional, $P_i$ has a weight decomposition $
P_i=C_i\oplus \bigoplus_n P_i^n. $ Then,  the components of
weight $n$, $$P^n:=\bigoplus_i P_i^n\ ,$$ form a subcomplex of
$P$, and the same is true for $C:=\bigoplus C_i$. So we have a
weight decomposition of $P$ as a direct sum of complexes
$$
P=C\oplus \bigoplus_n P^n.
$$
Taking homology we obtain
$$
HP=HC\oplus\bigoplus_n HP^n.
$$
Obviously, this decomposition is exactly the weight
decomposition of $HP$. Purity of $f$ implies $HC=0$ and
$H(P^n)=H_n(P)$, for all $n\in \mathbb Z$. Hence the inclusion
$
\bigoplus_{n} P^n\rightarrow P$ is a weak equivalence.

Next, for every $n\in \mathbb Z$, the homology of the complex
$P^n$ is concentrated  in degree $n$. So there is a natural
way to define a weak equivalence between the complex $P^n$ and
its homology $H(P^n)$. Let  $\tau_{\ge n}P^n$ be the canonical
truncation in degree $n$ of $P^n$,
$$
\tau_{\ge n}P^n := Z_n P^n\oplus \bigoplus_{i>n} P^n_{i}\ .
$$
This is a subcomplex of $P^n$, and the inclusion  $\tau_{\ge
n}P^n\rightarrow P^n$ is a weak equivalence. Since $\tau_{\ge
n}P^n$ is non trivial only in degrees $\ge n$, and its
homology is concentrated in degree $n$,  the canonical
projection $\tau_{\ge n}P^n
\rightarrow H(P^n)$ is a chain map, which is a weak equivalence.

Define $T$ by
$$
TP: =  \bigoplus_{n}\tau_{\ge n}P^n    \ .
$$
Obviously $T$ is an additive functor, and
$TP=\bigoplus_nT(P^n)$. Moreover, $T$ is a subfunctor of the
identity functor of $\CC^w(\mathbf k)$, and the canonical
projection $TP\lra HP$ is a weak equivalence.

We prove now that $T$ is a symmetric monoidal subfunctor of
the identity. Let $P$ and $Q$ be pure complexes of weight $0$.
Since $T$ is additive, and $ \sum_{i+j=n}P^i\otimes
Q^j\subset(P\otimes Q)^n$, it suffices to show that
$T(P^i)\otimes T(Q^j)\subset T(P^i\otimes Q^j)$. By the
Leibnitz rule,  we have an inclusion in degree $i+j$:
$$
Z_iP^i\otimes Z_jQ^j  \subset Z_{i+j}(P^i\otimes Q^j).
$$
In the other degrees the inclusion is trivially true. Hence,
$T$ being stable by products, it is a symmetric monoidal
subfunctor  of the identity.

Finally, the projection on the homology $TP\rightarrow HP$ is
well defined and obviously compatible with the K\"{u}nneth
morphism, so the canonical projection $T\rightarrow H$ is a
monoidal natural transformation. Therefore $C^w_*(\mathbf k) $
is a formal symmetric monoidal category.
\end{proof}

\begin{remark} The formality theorem for the current complex,
\ref{curretsformalfunctor}, could be obtained as a corollary  from the formality of the
full subcategory of $\CC_*(\mathbb R)$ whose objects are the
double complexes that satisfy the $dd^c$-lemma.
\end{remark}

\begin{corollary}\label{pesosimplicaformal} Let $P$ be an
operad with homology of  finite  type. If $P$  has a pure
endomorphism (with respect to some weight function $w$), then
$P$ is a formal operad.
\end{corollary}

\begin{proof}If
$P_\infty\rightarrow P$ is a minimal model of $P$, then
$P_\infty$ is an operad of finite type by \ref{finitudob}.
From the lifting property (\cite{MSS}, 3.123), there exists an
induced pure endomorphism $f$ on $P_\infty$. Thus
$(P_\infty,f)$ is an operad of $\CC^w_*(\mk)$, and the
corollary follows from th. \ref{categoriadepesosformal} and
prop.
\ref{formalfunctorscreateformaloperads}.
\end{proof}

Let  $P$ be an operad.  Since $\Sigma$-actions and
compositions $\circ_i$ are homogeneous maps of degree $0$,
every grading automorphism, with respect to a non root of unit
$\alpha$, is a pure endomorphism of the operad $HP$.

\begin{theorem} \label{aixeca}
Let $\mathbf k$ a field of characteristic zero, and  $P$ an
operad with homology of finite  type. The following statements
are equivalent:
\begin{enumerate}
\item $P$ is formal. \item There exists a model  $P'$ of $P$  such
that  $H: \Aut(P') \longrightarrow \Aut (HP)$ is surjective.
\item There exists a model $P'$ of $P$ and  $f\in Aut(P')$ such that $H(f)=\phi_{\alpha}$,
for some  $\alpha \in \mathbf k^*$  non root of unity.
 \item  There exists a pure
endomorphism $f$ in a model $P'$ of $P$.
\end{enumerate}
\end{theorem}

\begin{proof}
$(1)\Rightarrow (2)$ is prop. \ref{ooth}. $(2) \Rightarrow
(3)$ and $(3) \Rightarrow (4)$ are obvious. Finally,
$(4)\Rightarrow (1)$ is cor. \ref{pesosimplicaformal}.
\end{proof}

\section{Descent of formality}

In this section $\mathbf{k}$ will denote a field of
characteristic zero, and   operad will means an operad in  the
category $\CC_*(\mathbf k)$, unless another category was
mentioned. Using the characterization of formality of th.
\ref{aixeca}, we will prove now that formality does not depend
on the ground field, if it has zero characteristic.

\subsection{Automorphism group of a finite type operad}

\subsubsection{}
  Let $P$ be an operad. Restricting the
automorphism we have an inverse system of groups $
  \left(\Aut (t_nP) \right)_n
$ and a morphism of groups $\Aut (P) \longrightarrow \invlim
\Aut (t_nP)$. Because $P\cong
\lim\limits_{\leftarrow}t_{*}t_nP$, the following lemma is
clear.

\begin{lemma}\label{limitauts} The morphisms of restriction
 induce a canonical isomorphism of groups $$\Aut (P) \lra \invlim
\Aut (t_nP).$$
\end{lemma}

\subsubsection{}In order to prove that the group of automorphisms of a finite
type operad is an algebraic group,  we start by fixing some notations about group
schemes. Let  $\mathbf k\lra R$ be a commutative $\mathbf k$-algebra. If $P$ is an
operad, its extension of scalars $P\otimes_{\mathbf k} R$ is an operad in $\CC_*(R)$, and
the correspondence
$$
R \mapsto \AUT{P}(R) = \Aut_R (P \otimes_{\mathbf{k}} R) \ ,
$$
where $\Aut_R $ means the set of automorphisms of operads in
$\CC_*(R)$, defines a functor
$$
\AUT{P} : \mathbf{k}-\mathbf{alg} \longrightarrow \mathbf{Gr},
$$
from the category  $\mathbf{k}-\mathbf{alg} $  of commutative
$\mathbf k$-algebras, to the category
 $\mathbf{Gr}$ of groups. It is clear that
$$
\AUT{P} (\mathbf{k}) = \Aut (P) \ .
$$

We will denote by $\mathbb{G}_m$ the {\it multiplicative group
scheme\/} defined over the ground field $\mathbf{k}$.

\begin{proposition}
\label{algebraicos} Let  $P$ be a  truncated operad. If $P$ is of finite type, then
\begin{enumerate}
\item
$\Aut(P)$ is an algebraic matrix group over $\mathbf k$.
\item

$\AUT{P}$ is an algebraic affine group scheme over $\mathbf
k$, represented by the algebraic matrix group $\Aut (P)$.
\item Homology defines a morphism  $\mathbf{H} :
\AUT{P} \longrightarrow \AUT{HP}$  of algebraic affine group schemes.
\end{enumerate}
\end{proposition}

\begin{proof}  Let $P$  a finite type $n$-truncated operad.
The sum
 $M = \sum_{l \leq n} \mathrm{dim}\, P(l)$ is finite, hence
 $\Aut
(P)$ is the closed subgroup of $\GL{M}{k}$ defined by the
polynomial equations that express the compatibility with the
$\Sigma$-action, the differential, and the bilinear
compositions $\circ_i$. Thus $\Aut (P)$ is an algebraic matrix
group. Moreover, $\AUT{P}$ is obviously the algebraic affine
group scheme represented by the matrix group $\Aut(P)$.

Next, for every commutative $\mathbf k$-algebra $R$, the map
$$
\AUT{P}(R) = \Aut_R (P \otimes_{\mathbf{k}}R) \longrightarrow
 \Aut_R (HP \otimes_{\mathbf{k}}R) = \AUT{HP} (R)
$$
is a morphism of groups and it is natural in $R$; thus (3)
follows.
\end{proof}

\begin{theorem}\label{grupsalgebraics} Let $\mathbf k$ be a field of characteristic zero,
and $P$  a finite type truncated  operad. If $P$ is minimal,
then
$$
\mathbf{N} = \ker \left( \mathbf{H} :\AUT{P} \longrightarrow \AUT{HP} \right)
$$
is a unipotent algebraic affine group scheme over $\mathbf k$
.
\end{theorem}

\begin{proof} Since $\mathbf{k}$ has  zero characteristic, and $\AUT{P}$, $\AUT{HP}$ are
algebraic by prop. \ref{algebraicos}, $\mathbf{N}$ is
represented by an algebraic matrix group defined over
$\mathbf{k}$ (see \cite{Bo}). So it suffices to verify that
all elements in $\mathbf{N}(\mathbf{k})$ are unipotent.

Given $f \in \mathbf{N}(\mathbf{k})$, let $P^1=\ker
(f-\mathrm{id})^{\infty}$ be the primary component of $P$
corresponding to the eigenvalue $1$ (see \ref{weights}). Then
$P^1$ is a suboperad of $P$, and the inclusion $P^1
\hookrightarrow P$ is a weak equivalence. Since $P$ is
minimal, it follows from prop. \ref{almendruco} that $P =
P^1$, thus $f$ is unipotent.
\end{proof}

\subsection{A descent theorem}
After  these preliminaries, let us prove the descent theorem
of the formality for operads. In rational homotopy theory,
this corresponds to the descent theorem of formality for cdg
algebras of Sullivan and Halperin-Stasheff (\cite{Su} and
\cite{HS}, see also \cite{Mor} and \cite{R})

\begin{theorem}
\label{descens} Let $\mathbf{k}$ be a field of characteristic zero, and $\mathbf{k}
\subset \mathbf{K}$
  a field
extension. If $P$ is an operad in $\CC_*(\mathbf k)$ with
homology of finite type, then the following statements are
equivalent:
\begin{enumerate}
\item $P$ is formal.
\item $P\otimes \mathbf{K}$ is a formal operad in $\CC_*(\mathbf K)$.
\item For every $n$, $t_nP$ is formal.
\end{enumerate}
\end{theorem}

\begin{proof}
Because the statements of the theorem only depend on the
homotopy type of the operad, we can assume $P$ to be minimal
and, by th. \ref{finitudob}, of finite type.  Moreover,
minimality of $P$ is equivalent to the minimality of all its
truncations, $t_nP$.

Let us consider the following additional statement:

 $(2\frac{1}{2})$ \quad For every $n$, $t_nP\otimes\mathbf{K}$ is formal.

We will prove the following sequence of implications
$$
(1)\Rightarrow (2)\Rightarrow (2\tfrac{1}{2})\Rightarrow
(3)\Rightarrow (1),
$$

$(1)$ implies $(2)$  because $\_
\otimes_{\mathbf{k}}\mathbf{K} $ is an exact functor.

 If $P\otimes \mathbf{K}$ is formal, then so are all of
its truncations $t_n (P\otimes \mathbf K) \cong t_nP \otimes
\mathbf{K}$, because truncation functors are exact, so  $(2)$
implies $(2\frac{1}{2})$.

Let us see that $(2\frac{1}{2})$ implies $(3)$. From the
implication $(1)\Rightarrow (2)$, already proven, it is clear
that we may assume  $\mathbf{K}$ to be algebraically closed.
So, let $\mathbf{K}$ be an algebraically closed field, $n$ an
integer, and $P$ a finite type minimal operad such that  $t_nP
\otimes \mathbf{K}$ is formal. Since
$$
\AUT{t_nP}(\mathbf{K}) \longrightarrow \AUT{Ht_nP}(\mathbf{K})
$$
is a surjective map, by th. \ref{aixeca}, it results that
$$
\AUT{t_nP} \longrightarrow \AUT{Ht_nP}
$$
is a quotient map. Thus, by (\cite{Wa}, 18.1),  we have an
exact sequence of groups
$$
1 \longrightarrow \mathbf{N}(\mathbf{k}) \longrightarrow
\AUT{t_nP}(\mathbf{k})
\longrightarrow \AUT{Ht_nP} (\mathbf{k}) \longrightarrow H^1(\mathbf{K}/\mathbf{k},
\mathbf{N}) \longrightarrow \dots
$$
Since $\mathbf{N}$ is unipotent by th. \ref{grupsalgebraics},
and $\mathbf k$ has zero characteristic, it follows that
$H^1(\mathbf{K}/\mathbf{k}, \mathbf{N})$ is trivial ({\it op.~
cit.\/}, 18.2.e). So we have an exact sequence of groups
$$
1 \rightarrow \mathbf{N}(\mathbf{k}) \longrightarrow \Aut
(t_nP) \longrightarrow \Aut (Ht_nP ) \longrightarrow 1
$$
In particular, $\Aut (t_nP) \longrightarrow \Aut (Ht_nP)$ is
surjective. Hence, again by th. \ref{aixeca}, $t_nP$ is a
formal operad.

Let us see finally that $(3)$ implies $(1)$. By th.
\ref{aixeca}, it suffices to prove that all the grading
automorphisms have a lift. Let $\phi:\mathbb
G_m\lra\AUT{Ht_nP} $ the grading representation that sends
$\alpha\in \mathbb G_m$ to the grading automorphism
$\phi_\alpha$ defined in \ref{weights}. For every $n$, form
the pull-back of algebraic affine group schemes:
$$
\begin{CD}
\eg{n}       @>>>      \mathbb{G}_m           \\
@VVV                               @VV\phi V     \\
\AUT{t_nP} @>\mathbf H>>   \AUT{Ht_nP}
\end{CD}
$$
That is to say, for every commutative $\mathbf k$-algebra $R$,
$$
\eg{n}(R) = \mathrm{}\left\{ (f,\alpha ) \in \AUT{t_nP} (R) \times \mathbb{G}_m (R) \ ; \
Hf = \phi_{\alpha} \right\}\ .
$$
By \ref{limitauts}, we have a commutative diagram
$$
\begin{CD}
\lim\limits_\leftarrow\mathbf F_n(\mathbf k)@>>>\mathbb
G_m(\mathbf k)\\
@VVV@VV\phi V\\
\Aut(P)\cong\lim\limits_\leftarrow \Aut(t_nP)@>H>>\Aut(HP)\cong\lim\limits_\leftarrow
\Aut(Ht_nP)
\end{CD}
$$
so, to lift grading automorphisms, it suffices to verify that
the map $\lim\limits_\leftarrow\mathbf F_n(\mathbf
k)\longrightarrow\mathbb G_m(\mathbf k)$ is surjective. In
order to prove this surjectivity, first we will replace the
inverse system $\left(\mathbf F_n\right(\mathbf k))_n$ by an
inverse system $\left(\mathbf F'_n\right(\mathbf k))_n$ whose
transition maps are surjective. Indeed, for all $p \geq n$,
the restriction $\mathbf{\varrho}_{p,n} : \eg{p}
\longrightarrow \eg{n}$ is a morphism of algebraic affine
group schemes which are represented by   algebraic matrix
groups, so, by (\cite{Wa}, 15.1), it factors  as a quotient
map and a closed embedding:
$$
\eg{p} \longrightarrow \mathrm{im}\mathbf{\varrho}_{p,n} \longrightarrow \eg{n} \ .
$$
Denote $\egprima{n} :=\bigcap_{p \geq n}
\mathrm{im}\mathbf{\varrho}_{p,n}$.
 Since $\{
\mathrm{im}\mathbf{\varrho}_{p,n} \}_{p \geq n}$ is a descending chain of closed
subschemes of the noetherian scheme $\eg{n}$, there exists an
integer $N(n) \geq n$ such that
$$
\egprima{n} = \mathrm{im} \mathbf{\varrho}_{N(n),n} \ ,
$$
thus the restrictions $\mathbf{\varrho}_{n+1,n}$ induce quotient maps
$\mathbf{\varrho}_{n+1,n} : \egprima{n+1} \longrightarrow \egprima{n}$. So, applying
again (\cite{Wa}, 18.1),  we have an exact sequence of groups
$$
1 \longrightarrow \mathbf{N}'(\mathbf{k}) \longrightarrow
\egprima{n+1}(\mathbf{k})
\longrightarrow \egprima{n}(\mathbf{k}) \longrightarrow H^1(\clausura /\mathbf{k},
\mathbf{N}) \longrightarrow \dots
$$
Here, $\mathbf{N}'(\mathbf{k})$ is a closed subscheme
 of $\mathbf{N}(\mathbf{k})$
because, for every $(f, \alpha) \in \mathbf{N}'(\mathbf{k})$
we have  $\alpha = 1$ and so $Hf = 1$ in $Ht^*_{\le n+1}P$,
which means that $f \in \mathbf{N}(\mathbf{k})$. By th.
\ref{grupsalgebraics}, $\mathbf{N}'(\mathbf{k})$ is unipotent, thus, as in the previous
implication, it follows that $\egprima{n+1}(\mathbf{k})
\longrightarrow \egprima{n} (\mathbf{k})$ is surjective for
all $n \geq 2$.

Since in the inverse system $\left(\mathbf F'_n\right(\mathbf
k))_n$ all the transition maps are surjective, the map
$$
\invlim \egprima{p}(\mathbf{k} )\longrightarrow \egprima{2} (\mathbf{k} ) \ .
$$
is surjective as well.  Moreover, $\egprima{2}(\mathbf{k})
\longrightarrow
\mathbb{G}_m(\mathbf k)$ is also surjective. Indeed, given $\alpha \in
\mathbb{G}_m(\mathbf k)$, since $t_{ N(2)}{P}$ is formal by hypothesis,  by th.
\ref{aixeca}  we can lift the grading automorphism $\phi_{\alpha} \in \Aut (Ht_{ N(2)}P)$
to an automorphism $f\in \Aut (t_{ N(2)}{P})$. So we have an
element $(f, \alpha) \in
\eg{N(2)}(\mathbf{k})$, whose image in $\eg{2}(\mathbf{k})$ will be an element of
$\egprima{2} ( \mathbf{k})$ which will project onto $\alpha$.

We conclude that $\invlim \egprima{p}(\mathbf{k} )
\longrightarrow \mathbb{G}_m(\mathbf k) $ is surjective, hence
$P$ is formal.
\end{proof}

\subsection{Applications} As an immediate consequence of th.
\ref{descens}, the previous  theorems \ref{cadenesformalsR}, and \ref{Coperadformal}  of
formality over $\mathbb R$ imply, respectively,
 the
following corollaries

\begin{corollary}\label{cadenesformalsQ}
If $X$ is an operad in $\CKM$, then the operad of cubic chains
$C_*({X}; \mathbb Q)$ is formal.
\end{corollary}

\begin{corollary} If $X$ is an
operad in $\Vman$, then the operad of cubic chains
$C_*(X;\mathbb{Q})$ is formal.
\end{corollary}

Finally, we can apply th. \ref{descens} to the formality of
the little $k$-disc operad. Let $\mathcal D_k$  denote the
{\it little $k$-discs \/} operad of Boardman and Vogt. It is
the topological operad with $\mathcal D_k (1) = \mathbf{pt} $,
and, for  $l\geq 2$, $\mathcal D_k (l) $ is  the space of
configurations of $l$ disjoint discs inside the unity disc of
$\mathbb{R}^k$.

M. Kontsevich  proved that the operad  of cubic chains
$C_*(\mathcal D_k,\mathbb R) $ is formal (\cite{Ko}).
Therefore, from th. \ref{descens}, we obtain

\begin{corollary} The operad of  cubic chains of the little
$k$-discs operad $C_*(\mathcal D_k ;\mathbb{Q})$ is  formal.
\end{corollary}

\section{Cyclic operads}

\subsection{Basic results} Let us recall some definitions from \cite{GeK94} (see also
\cite{GeK98} and \cite{MSS}). For all $l\in \mathbb N$, the group
$\Sigma^+_l:=\Aut\{0,1,\dots,l\}$ contains $\Sigma_l$ as a
subgroup, and it is generated by $\Sigma_l$ and the cyclic
permutation of order $l+1$,\quad
$\tau_l:(0,1,\dots,l)\mapsto(1,2,\dots,l,0)$.

Let $\mathcal C$ be a symmetric monoidal category. A {\it
cyclic\/} $\Sigma$-module $E$ in  $\mathcal{C}$ is a sequence
$(E(l))_{l\ge 1}$ of objects of $\mathcal C$ together with an
action of $\Sigma^+_l $ on each $E(l)$. Let $\sModplus$ denote
the category of cyclic $\Sigma$-modules. Forgetting the action
of the cyclic permutation we have a functor
$$
U^-: \sModplus \longrightarrow \sMod \, .
$$

A {\it cyclic operad\/} is a cyclic $\Sigma$-module $P$ whose
underlying $\Sigma$-module $U^-P$ has the structure of an
operad compatible with the  action of the cyclic permutation
(see {\it loc. cit.}). Let $\Opplus$ denote the category of
cyclic operads. We also have an obvious forgetful functor
$$
U^- : \Opplus \longrightarrow \Op \ .
$$

There are obvious extensions of the notions of free operad, homology, weak equivalence,
 minimality and formality for cyclic dg operads and all the results
in the previous sections can be easily transferred to the cyclic setting. In particular,
every cyclic dg operad $P$ with $HP(1)=0$,  has a minimal model $P_\infty$. Moreover
$U^-(P_\infty)$ is a minimal model of $U^-(P)$. Finally, we can deduce results analogous
to the formality criterion (th. \ref{aixeca}), and to  the descent of formality (th.
\ref{descens}) for cyclic operads.

Let $\mathcal A$ an abelian category. It is clear that a formal symmetric monoidal
functor $F: \mathcal{C} \longrightarrow \Chainsab$ induces a functor of cyclic operads
$$
F: \Opplus\!_{\mathcal{C}} \longrightarrow
\Opplus\!_{\Chainsab}
$$
which sends cyclic operads in $\mathcal{C}$ to formal cyclic
operads in $\Chainsab$. From theorems \ref{cadenesformals} and
\ref{descens} it follows that $C_*(X; \mathbb{Q})$ is a formal
cyclic operad, for every cyclic operad $X$ in $\CKM$.

\subsection{Formality of the cyclic operad $C_*(\overline{\M}_{0};\mathbb Q)$ }
Let us apply the  previous  results  to the { configuration
operad\/}.

Let $\M_{0,l}$ be the moduli space  of $l$ different labelled
points on the complex projective line ${\mathbb P}^1$. For $l
\geq 3$, let $\overline{\M}_{0,l}$ denote its {\it
Grothendieck-Knudsen compactification\/}, that is, the moduli
space of stable curves of genus $0$, with $l$ different
labelled points.

For $l=1$, put $\modulizero (1) = *$, a point, and for $l \geq
2$, let $\modulizero (l) =
\overline{\M}_{0,l+1}$. The family of spaces $\modulizero = \left(\modulizero
(l)\right)_{l\geq 1}$ is a cyclic  operad in $\Var$
(\cite{GK},  or \cite{MSS}). Applying the functor of cubic
chains componentwise we obtain a dg cyclic operad
$C_*(\modulizero ;
\mathbb{Q} )$. So, we have the following results.

\begin{corollary}\label{modulispaceformal} The cyclic operad of cubic
chains  $C_*(\modulizero ; \mathbb{Q})$  is formal.
\end{corollary}

\medskip

\begin{corollary} The categories of strongly homotopy
$C_*(\modulizero ; \mathbb{Q})$ and $H_*(\modulizero ;
\mathbb{Q})$-algebras are equivalent.
\end{corollary}

\section{Modular operads}

\subsection{Preliminaries}
Let us recall some definitions and notations about modular
operads (see \cite{GeK98}, or
\cite{MSS}, for details).

\subsubsection{} Let $\mathcal C$ be a symmetric monoidal
category. A {\it modular $\Sigma$-module\/} of $\mathcal C$ is a bigraded object of
$\mathcal{C}$, $E=\left(E((g,l))\right)_{g,l}$, with $g,l\ge0,\, 2g-2+l>0$, such that
$E((g,l))$ has a right $ \Sigma_{l}$-action. Let us denote the category of modular
$\Sigma$-modules by $\MMod_{\mathcal{C}}$, or just $\MMod$ if no confusion can arise.

\subsubsection{}
A {\it modular operad\/} is a  modular $\Sigma$-module $ P$, together with {\it
composition} morphisms
$$
\circ_i : P((g,l)) \otimes P((h,m)) \longrightarrow P((g+h,l+m-2)) \ ,\; 1\le i\le l,
$$
and    {\it contraction\/} morphisms
$$
\xi_{ij}: P((g,l)) \longrightarrow P((g+1,l-2)) \ , \;1\le i\not= j\le l,
$$
which verify axioms of associativity, commutativity and compatibility (see \cite{GeK98},
\cite{MSS}). Let us denote the category of modular operads by $\MOp_{\mathcal{C}}$, or
just $\MOp$ if no confusion can arise.

 As for operads and cyclic operads, from the definitions it follows  that every symmetric
monoidal functor $F : \mathcal{C} \longrightarrow \mathcal{D}$
applied componentwise induces a functor
$$
\MOp_{F}:\MOp_{\mathcal{C}} \longrightarrow \MOp_{\mathcal{D}} \ ,
$$
and every monoidal natural transformation $\phi:F\Rightarrow
G$ between symmetric monoidal functors induces a natural
transformation $\MOp_\phi:\MOp_{F}\Rightarrow
\MOp_G.$

\begin{example}
As Getzler and Kapranov proved (\cite{GeK98}),  the family
$\overline{\M}((g,l)):=\overline{\M}_{g,l}$ of
Deligne-Knudsen-Mumford moduli spaces  of stable genus $g$
algebraic curves with $l$ marked points, with the maps that
identify marked points, is a modular operad in the category of
projective  smooth DM-stacks.
\end{example}
 \subsection{dg modular operads} \label{functorformalmodular}

{\it From now on, $\mathbf k$ will denote a field of
characteristic zero, and modular operads in $\CC_*(\mathbf k)$
will be called simply dg modular operads.} \vskip 3mm

\begin{example}
Let $V$ be a finite type chain complex of $\mathbf k$-vector
spaces, and $B$ an inner product  over $V$, that is to say, a
non-degenerate graded symmetric bilinear form $B:V\otimes
V\lra \mathbf k$ of degree $0$.  It is shown in \cite{GeK98}
that there exists a dg modular operad $\mathcal E[V]$ such
that
$$\mathcal E[V]((g,l))=V^{\otimes l},$$
with the obvious structure morphisms.
\end{example}

An {\it ideal\/} of a dg modular operad $P$ is a modular $\Sigma$-submodule $I$ of $P$,
such that  $P\cdot I\subset I$,  $I\cdot P\subset I$,  and $I$ is closed under the
contractions $\xi_{ij}$.

For any dg modular operad $P$ and any ideal $I $ of $P$, the
quotient $P/I$, inherits a natural structure of dg modular
operad and the projection $P \longrightarrow P/I$ is a
morphism of dg modular operads.

If $P$ is a dg modular operad its {\it homology} $HP$, defined by
$(HP)((g,l))=H(P((g,l)))$, is also a dg modular operad. A morphism $\rho:P\lra Q$ of dg
modular operads is said to be a {\it weak equivalence} if $\rho((g,l)):P((g,l))\lra
Q((g,l))$ is a weak equivalence for all $(g,l)$.

The localization of $\MOp_{\CC_*(\mathbf k)}$ with respect to
the weak equivalences is denoted by $\mathrm{ Ho}
\MOp_{\CC_*(\mathbf k)}$.

\begin{definition}
A dg modular operad $P$ is said to be {\it formal} if $P$ is
weakly equivalent to its homology $HP$.
\end{definition}

 Clearly, for a formal symmetric
monoidal functor $F : \mathcal{C} \longrightarrow
\CC_*(\mathbf k)$, the induced functor
$$
F: \MOp_{\mathcal{C}} \longrightarrow \MOp_{\CC_*(\mathbf k)}
$$
transforms modular operads in $\mathcal{C}$ to  formal modular
operads in $\CC_*(\mathbf k)$.

\subsection{Modular dimension}
In order to study the homotopy properties of dg modular operads we will replace the arity
truncation
 with the truncation with respect to the modular dimension.

 Let $\mathcal C$ be a symmetric monoidal category.

Recall that, the  dimension as algebraic variety of the moduli
space $\overline{\M}_{g,l}$   is $3g-3+l$. So the following
definition is a natural one. The function $d:\mathbb Z^2\lra
\mathbb Z$, given by $d(g,l)=3g-3+l$, will  be called the {\it
modular dimension} function.

\begin{definition} Let $E$ be
a modular $\Sigma$-module in $\mathcal C$. The modular dimension function induces a
graduation $\left( E_n\right)_{n\ge 0}$ on $E$ by $E_n=\left(E((g,l))\right)_{d(g,l)=n}$,
and a decreasing filtration $\left(E_{\ge n}\right)_n$ of $E$ by
$$
E_{\ge n}:=\left(E((g,l))\right)_{d(g,l)\ge n} \ .
$$
\end{definition}

 The following properties are easily checked.

\begin{proposition}\label{propiedaddimensionmodular}
 Let $P$ be a modular operad in $\mathcal C$. The modular dimension grading
satisfies
$$P\cdot P_n\subset P_{\ge n+1},$$ where
$P\cdot P_n$ is the set of meaningful products
$\alpha\circ_i\beta$, with $\alpha\in P_m$, $\beta\in P_l$ and
at least one of \,$l,\,m$\, is \,$n$. On the other hand, the
contraction maps satisfies
$$
\xi_{ij}: P_n \longrightarrow P_{\ge n+1} \ ,
$$
for all $i,j$.
\end{proposition}

\subsection{Truncation of modular operads}

\begin{definition}

 A {\it $n$-truncated modular operad} in a symmetric monoidal
category $\mathcal C$ is a modular operad defined only up to
modular dimension $n$, that is, a family of $\mathcal{C}$,
$\{P((g,l));\,g,l\ge 0, \; 2g-2+l>0,\, d(g,l)\le n\,\}$,  such
that $P((g,l))$ has a right $
\Sigma_{l}$-action, with morphisms
$$
\circ_i : P((g,l)) \otimes P((h,m)) \longrightarrow P((g+h,l+m-2)) \ ,\; 1\le i\le l,
$$
and    {\it contractions\/}
$$
\xi_{ij}: P((g,l)) \longrightarrow P((g+1,l-2)) \ , \;1\le i\not= j\le l,
$$
satisfying those axioms  in $\MOp$ that make sense.

\end{definition}

If $\MOp_{\le n}$ denotes the category of $n$-{\it truncated
dg modular operads\/}, we have a {\it modular dimension
truncation\/} functor,
$$
t_{ n}:\MOp \longrightarrow \MOp_{\le n}\,,
$$
defined by $t_n(P)=(P((g,l)))_{d(g,l)\le n}$.

Since the obvious forgetful functor
$$
U: \MOp \longrightarrow \MMod \
$$
has a left adjoint, the {\it free modular operad functor\/}
$$
\mathbb{M} : \MMod \longrightarrow \MOp \ ,
$$
(see\cite{GeK98}, 2.18), by prop. \ref{propiedaddimensionmodular} we can translate the
truncation formalism developed in \ref{Truncatedoperads} to the setting of dg modular
operads. So we have a sequence of adjunctions $t_{!}\dashv t_{n}\dashv t_{*}$, and the
propositions \ref{propiedadest*} and \ref{propiedadest!} are still true, merely replacing
``operad" with ``modular operad", and ``arity" with ``modular dimension" shifted by $+2$.
For instance, the arity truncation begins with $t_2$, whereas the modular dimension
truncation begins with $t_0$.

If $P$ is a dg modular operad, the direct system of dg modular operads given by
$$
0\lra t_{!}t_{0}P\rightarrow\cdots \rightarrow t_{!}t_{n-1} P
\rightarrow t_{!}t_{n}P\rightarrow\cdots
$$
is called the {\it canonical tower} of $P$.

\subsection{Principal extensions}

Let us  explicitly describe the construction of a principal extension in the context of
modular operads. Let $P$ be a dg modular operad,  $V$  a dg modular $\Sigma$-module with
zero differential, concentrated in modular dimension $n\ge 0$, and $\xi:V[-1]\lra P_n$ a
chain map. Then the principal extension of $P$ by $\xi$, $P\sqcup_\xi V$,  is defined by
a universal property as in \ref{estructuradeunaextensionprincipal},
$$
\Hom_{\MOp}(P\sqcup_\xi V,Q)=  \{(f,g);\ f\in \Hom_{\MOp}(P, Q)\ ,\, g\in \Hom_{\bf
GrMMod}(V, UQ)_{0}\ , \ d_Qg-gd_V=f\xi\ \}\, .
$$

In particular we have
 $$(P\sqcup_\xi
V)_i=\left\{
\aligned \label{}&P_i,\hskip 3mm \text{ if } i<n,\\
 &P_n\oplus_\xi V,\hskip 3mm \text{ if } i=n,\endaligned
\right.$$
 because in $t_{ n}(P\sqcup_\xi
V)$ all  the structural morphisms
     involving $V$ are trivial, by prop.
     \ref{propiedaddimensionmodular}.

Furthermore, the following property, analogous  to prop.
\ref{estructuradeunaextensionprincipal}, is satisfied.
\begin{proposition}
\label{estructuradeunaextensionprincipaldemodularoperads} Let $n\ge 0$ an integer. Let
$P$ be a dg modular operad such that $t_!t_nP\cong P$, $V$ a dg modular $\Sigma$-module
concentrated in modular dimension $n$, with zero differential, and $\xi:V[-1]\lra P_n$ a
morphism of dg modular $\Sigma_n$-modules. The principal extension $P\sqcup_\xi V$
satisfy

(1)\quad  $t_{n-1}(P\sqcup_\xi V)\cong t_{n-1}P$.

(2)\quad  $(P\sqcup_\xi V)_n\cong C\xi,$ in particular, there
exists an exact sequence of complexes
$$
0\lra P_n\lra (P\sqcup_\xi V)_n\lra V\lra 0.
$$
(3) \quad $P\sqcup_\xi V\cong t_{!}t_n(P\sqcup_\xi V) $.

(4) A morphism of dg modular operads $\phi:P\sqcup_\xi V\lra
Q$ is determined by a morphism of $n$-truncated dg modular
operads $f:t_nP\lra t_nQ$, and a homogeneous map $g:V\lra Q_n$
of modular $\Sigma$-modules, such that $f\xi =dg$.
\end{proposition}

\subsection{Minimal models}

\subsubsection{Minimal objects}
\begin{definition} A dg modular operad $M$ is said to be {\it minimal\/}
if  the canonical tower
$$
0\longrightarrow t_{!}t_{0}M\longrightarrow \dots
\longrightarrow t_{!}t_{ n-1}M
\longrightarrow t_!t_{ n}M \longrightarrow \dots
$$
is a sequence of principal extensions.

 A {\it minimal model\/} of a dg modular operad $P$ is a minimal dg modular operad
$P_\infty$ together with a weak equivalence $P_\infty
\longrightarrow P$ in $\MOp$.
\end{definition}

From prop. \ref{estructuradeunaextensionprincipaldemodularoperads}, it follows,  by
induction on the modular dimension, that:
\begin{proposition}
\label{quis=isomodular} Let $M$, $N$ be minimal dg modular operads. If $\rho:M
\longrightarrow N$ is a weak equivalence of dg modular operads, then $\rho$ is an
isomorphism.
\end{proposition}

\subsubsection{Existence of minimal models}

\begin{theorem}\label{recetapaella}
Let $\mathbf k$ be a field of characteristic zero. Every
modular operad $P$ in $\CC_*(\mathbf k)$ has a minimal model.
\end{theorem}
\begin{proof}
 We start in modular dimension $0$. Let
$M^{0}=\mathbb{M}HP_0$, and $s:HP_0\lra ZP_0$ a
 section of the canonical projection. Then $s$ induces a morphism
of modular operads
$$
\rho^{0}: M^{0} \longrightarrow P
$$
which is a weak equivalence of modular operads up to modular
dimension $0$, because $M^{0}_0=HP_0$.

For $n\ge 1$, assume that we have already constructed a
morphism of modular operads
$$
\rho^{n-1} : M^{n-1} \longrightarrow P
$$
such that
\begin{enumerate}
\item $M^{n-1}\cong t_{!}t_{n-1}M^{n-1}$ is a
minimal  modular operad, and
\item $t_{ n-1}(\rho^{n-1})$ is a weak
equivalence.
\end{enumerate}
To define the next step of the induction we will use the following statement, which
contains the main  homological part of the inductive construction of  minimal models.
\begin{lemma}
Let
$$
\begin{CD}
B@>\eta>>A@.@.\\
@V\lambda VV@VV\mu V @.@.\\
(C\zeta)[-1]@>-p_Y>>Y@>\zeta>>X
\end{CD}
$$
be a  commutative diagram  of complexes of an additive
category, then  there exists a chain map $\nu:C\eta\lra X$
such that in the
 diagram
$$
\begin{CD}
@.@. B @>\eta>> A @>>> C\eta @>>> B[1] \\
(8.6.3.1)\hskip 30mm @.@. @V\lambda VV
@VV\mu  V@VV\nu V @VV\lambda[1] V\\
@.@.(C\zeta)[-1]@>-p_Y>>Y @>\zeta>>X@>>>C\zeta
\end{CD}\hskip 50mm
$$
the central square is  commutative, and  the right hand side
square is homotopy commutative. Moreover, the rows of
$(8.6.3.1)$ are distinguished triangles, and the vertical maps
define a morphism of triangles in the derived category.
\end{lemma}

We have $C\eta=A\oplus_\eta B[1]$, and $C\zeta=X\oplus_\zeta
Y[1]$. Let $(\lambda_X,\lambda_Y)$ be the components of
$\lambda$, then one can check that $
\nu(a,b)=\lambda_X(b)+\zeta\mu(a), $  with the homotopy $
 h(a,b)=(0,\mu (a))$,
satisfies the conditions of the statement.

The upper row of the diagram $(8.6.3.1)$ is obviously a
distinguished triangle. By
 axiom
$(TR2)$ of a triangulated category, turning the distinguished
triangle
$$
\begin{CD}
Y@>\zeta>>X@>>>C\zeta@>p_Y>>Y[1]
\end{CD}
$$
one step to the left we obtain that the lower row of the diagram (8.6.3.1) is also a
distinguished triangle.

Now we return  to the proof of the theorem. Since $\mathbf k$
is a field of zero characteristic, the category of modular
$\Sigma$-modules is semisimple, and $C\rho^{n-1}_n$ is a
formal complex of modular $\Sigma$-modules. Therefore, if  $V=
HC\rho^{n-1}_n$, with the zero differential,  there exists a
weak equivalence
$$s:V\lra C\rho^{n-1}_n\,.
$$
In fact, $s$ can be obtained from a $\Sigma$-equivariant
section  of the canonical projection from cycles to homology.

Let $\xi$ be the composition
$$\begin{CD}
{V[-1]}@>s[-1]>> \left(C\rho^{n-1}_n\right)[-1]
@>-p>>M^{n-1}_n,
\end{CD}
$$
where the second arrow is the opposite of the canonical
projection.
 We have a commutative diagram of complexes
$$
\begin{CD}
V[-1]@>\xi>>M^{n-1}_{n}@.\\
@V s [-1]VV@VV \mathrm{id} V\\
\left(C\rho^{n-1}_n\right)[-1]@>-p>> M^{n-1}_n@>\rho^{n-1}_n>> P_n
.\\
\end{CD}
$$
By the previous lemma,  there exists a chain map
$$\nu:C\xi\lra P_n$$ that completes the previous diagram
in a diagram
$$
\begin{CD}
@.@. V[-1]@>\xi>>M^{n-1}_n@>>>C\xi@>>>
V@.@.\\
(8.6.3.2)\hskip 30mm@.@.@Vs[-1]VV@V\mathrm{id}VV@V\nu VV@VVsV@.@.\hskip 20mm
\\
@.@.\left(C\rho^{n-1}_n\right)[-1]@>-p>> M^{n-1}_n
@>\rho^{n-1}_n>> P_n@>>> C\rho^{n-1}_n.@.@.\\
\end{CD}\hskip 20mm
$$
where  the rows are distinguished triangles in the category of
complexes, the central square is commutative, and the vertical
maps define a morphism of triangles in the derived category.

The step $M^{n}$  is defined as the  principal extension of
$M^{n-1}$ by the attachment map $\xi:V[-1]\lra M^{n-1}_n$,
$$
M^{n}:=M^{n-1}\sqcup_{\xi} V\,.
$$
Let $\nu_V:V\lra P$ be the graded map
$$\begin{CD} V@>>> C\xi@>\nu>>P_n\end{CD}$$
where the first map is the canonical inclusion. Since
\,$\rho^{n-1}\xi=d \nu_{V}$\ ,  the maps $\rho^{n-1}$ and
$\nu_{V}$ define, according to the universal property of
$M^{n}=M^{n-1}\sqcup_{\xi} V$, a morphism of modular operads
$$
\rho^{n}:M^{n}\lra P
$$
such that $t_{ n-1}\rho^{n}=t_{ n-1}\rho^{n-1}$ and
$\rho_n^{n}=\nu$.  By the inductive hypothesis, $\rho^{n}$ is
a weak equivalence in modular dimensions $<n$. Finally, in the
diagram $(8.6.3.2)$  $s$ is a weak equivalence, hence $\nu$ is
a weak equivalence as well. It follows that $t_{n}\rho^{n}$ is
a weak equivalence, which finishes the induction. Therefore,
$\lim\limits_{\rightarrow}M^{n}$ is a minimal model of $P$.
\end{proof}

\subsubsection{Finiteness of minimal models}

\begin{definition}   A modular
$\Sigma$-module $V$  is said to be {\it of finite type\/} if, for every $(g,l)$,
$V((g,l))$ is a finite dimensional $\mathbf{k}$-vector space. A  dg modular operad $P$
is said to be of {\it finite type} if $UP$ is of finite type.
\end{definition}
Obviously, for every integer $n\ge 0$, there are only a finite
number of  pairs $(g,l)$ such that $g,l,2g-2+l>0$ and
$d(g,l)=n$, thus a modular $\Sigma$-module $V$ is of finite
type if, and only if, $V_n$ is finite dimensional, for every
$n\ge 0$.

\begin{proposition}\label{finitudmodular}
 If  $V$ is  a modular $\Sigma$-module of finite type, then $\mathbb M(V)$
is of finite type. \end{proposition}
\begin{proof}
Indeed,    for every  pair $(g,l)$, there is an isomorphism
$$
\mathbb M(V)((g,l))\cong \bigoplus_{\gamma\in \{\mathbf\Gamma((g,l))\}
}V((\gamma))_{\text{Aut}(\gamma)}
$$
where $\{\mathbf\Gamma((g,l))\} $ denotes the set of
equivalence classes of isomorphisms of stable $l$-labelled
graphs of genus $g$, the subscript $\text{Aut}(\gamma)$
denotes the space of coinvariants, and
$$
V((\gamma))=\bigotimes_{v\in
\text{Vert}(\gamma)}V((g(v),\text{Leg }{(v)})).
$$
By  \cite{GeK98} lemma 2.16, the set $\{\mathbf\Gamma((g,l))\}
$ is finite, for every  pair $(g,l)$. Therefore the free
modular operad $\mathbb M(V)$ is of finite type.
\end{proof}

As a consequence of prop. \ref{finitudmodular}  we obtain the
finiteness result analogous to \ref{finitudob}.
\begin{theorem}
Let $P$ be  a dg modular operad. If $HP$ is of finite type,
then every
 minimal model   of $P$  is of finite type.
\end{theorem}

\subsection{Lifting properties}
Analogous to def. II.3.121 of \cite{MSS}, there exists a similarly defined path object
and a notion of homotopy in the category of dg modular operads.

\subsubsection{Homotopy}
Let  $\mathbf I:=\mathbf k[t,\delta t]$ be the differential
graded commutative $\mathbf k$-algebra generated by a
generator $t$ in degree $0$ an its differential $\delta t$ in
degree $-1$. For every dg modular operad $P$, the path object
of $P$
 is   the dg modular operad
$P\otimes\mathbf I $, obtained by  extension of scalars.

 The evaluations at $0$ and $1$ define two morphisms of modular operads
 $\rho_0,\rho_1:P\otimes\mathbf I\rightrightarrows P$ which are
weak equivalences. An {\it elementary homotopy} between two morphisms of dg modular
operads $f_0,f_1:P\rightrightarrows Q$ is a morphism $H:P\lra Q\otimes\mathbf I$ of dg
modular operads such that $\rho_i H=f_i$, for $i=0,1$.  Elementary homotopy is a
reflexive and symmetric relation, and the {\it homotopy relation} between morphisms is
the equivalence relation generated by elementary homotopy. Homotopic morphisms induce the
same morphism in $\mathrm{ Ho }\MOp$.

\subsubsection{Lifting properties of minimal objects}
Obstruction theory, that is,  lemma II.3.139 {\it op. cit.}, and its consequences: the
homotopy properties of the minimal objects (theorems II.3.120 and II.3.123 {\it op.
cit.}), is easily established in the context of modular operads. So we have
\begin{lemma}\label{liftingproperty}
Let $\rho:Q\lra R$ be a weak equivalence of dg modular
operads, and $\iota : P
\longrightarrow P \sqcup_\xi V$ a principal extension. For every homotopy commutative
diagram in $\MOp$
$$
\begin{CD}
P                    @>\phi >>           Q       \\
@V\iota VV                            @V\rho VV   \\
P \sqcup_\xi V @>\psi >> R
\end{CD}
$$
there exists an extension $\overline{\phi} : P \sqcup_\xi V
\longrightarrow Q$ of $\phi$ such
 that $\rho \overline{\phi}$ is homotopic to $\psi$. Moreover, $\overline{\phi}$ is
 unique up to homotopy.
\end{lemma}

From this lemma  the lifting property of minimal modular
operads follows by induction:

\medskip

\begin{theorem}
\label{liftingpropertyforminimalsmodular} Let  $\rho : Q \longrightarrow R$ be a weak
equivalence of dg modular operads, and $M$ a minimal modular
operad. For every morphism $\psi : M \longrightarrow R$, there
exists a morphism $\widetilde{\psi} : M
\longrightarrow Q$ such that $\rho\widetilde{{\psi}}$ is homotopic to $\psi$. Moreover,
$\widetilde{\psi}$ is unique up to homotopy.
\end{theorem}

\subsubsection{Uniqueness of minimal models}

From th. \ref{liftingpropertyforminimalsmodular}  and prop. \ref{quis=isomodular}, we
obtain
\begin{theorem}\label{unicidadmmmodular}
Two minimal models of a modular operad are isomorphic.
\end{theorem}

The modular analogue of prop. \ref{almendruco} follows in the same way.

\begin{proposition}\label{almendrucomodular} Let $M$ be a minimal dg
modular operad and $P$ a subobject of $M$. If  the inclusion
$P \hookrightarrow  M$ is a weak equivalence, then $P = M$.
\end{proposition}

\subsection{Formality}

From theorems \ref{liftingpropertyforminimalsmodular} and \ref{quis=isomodular}, the
modular analogue of prop. \ref{ooth} follows easily.
\begin{proposition} \label{oothmodular} Let $M$ be a minimal dg modular operad. If $M$
is  formal, then the map $H:\nolinebreak\Aut(M)\longrightarrow\Aut(HM)$ is surjective.
\end{proposition}

Now, from th. \ref{pesosimplicaformal} and prop.
\ref{oothmodular}, the formality criterion for modular operads
follows with the same proof as th. \ref{aixeca}.

\begin{theorem} \label{aixecamodular} Let $\mathbf k$ be a field of characteristic zero,
and  $P$  a dg modular operad with homology of
   finite  type.
The following statements are equivalent:
\begin{enumerate}
\item $P$ is formal.
\item
There exists  a model  $P'$ of $P$  such that  $H: \Aut(P')
\longrightarrow \Aut (HP)$ is surjective.
\item There exists a model  $P'$ of $P$, and  $f\in \Aut(P')$ such that
$Hf=\phi_{\alpha}$, for some  $\alpha\in \mathbf{k^*} $ non
root of unity.
\item  There exists a pure endomorphism $f$ in  a  model
$P'$  of $P$.
\end{enumerate}
\end{theorem}
Then, using this result, the descent of formality for modular
operads follows as th. \ref{descens}.

\begin{theorem}\label{descensmodular} Let $\mathbf{k}$ be a field of
characteristic zero, and  $\mathbf{k} \subset \mathbf{K}$
  a field
extension. If $P$ is a  modular operad in $\CC_*(\mathbf k)$
with homology of finite type, then $P$ is formal if, and only
if, $P\otimes \mathbf{K}$ is a formal modular operad in
$\CC_*(\mathbf K)$.
\end{theorem}

Finally, the  result below  follows from \ref{functorformalmodular}, and theorems
\ref{hodgeforVman}, \ref{descensmodular}.
\begin{theorem}\label{Coperadmodformal} Let $X$ be a modular
operad in $\Vman$. Then $C_*(X; \mathbb{Q})$ is a formal
modular operad.
\end{theorem}

\subsection{Strongly homotopy algebras over a modular operad}
Let  $P$ be  a dg modular operad. Recall (\cite{GeK98}) that a
$P$-algebra is a finite type chain complex $V$ with an inner
product $B$,
 together with  a morphism of modular operads $P
\longrightarrow \mathcal E[V]$.

We give the following definition. A {\it strongly homotopy
$P$-algebra}, or {\it sh $P$-algebra}, is a finite type chain
complex $V$ with an inner product $B$,
 together with  a  morphism  $P
\longrightarrow \mathcal E[V]$ in  $\mathrm {Ho}\MOp$. By
\ref{liftingpropertyforminimalsmodular}, this is equivalent to giving a homotopy class of
morphisms $P_\infty\lra  \mathcal E[V]$.

Let $(V,B)$, and  $(W,B)$ be sh $P$-algebras. A {\it morphism
of sh $P$-algebras} $f$ is a chain map  $f:V\lra W$ compatible
with the inner products and
 such that the following diagram
$$\begin{array}{ccc}
   &  & \mathcal E[V] \\
   & \nearrow & \\
  P &  & \downarrow \; f_*  \\
   & \searrow & \\
   &  & \mathcal E[W]
\end{array}$$
 commutes in $\mathrm{Ho }\MOp$.
\medskip

The homotopical invariance  is an immediate consequence  of
the above definitions:
\begin{proposition} Let $(V,B)$ be a finite type
chain complex with an inner product, $(W,B)$ a sh $P$-algebra,
and $f: (V,B)
\longrightarrow (W,B)$ a chain map compatible with $B$ and such that $f:V\lra W$ is a
homotopy equivalence. Then $(V,B)$ has a unique structure of
sh $P$-algebra such that $f$ becomes a morphism of sh
$P$-algebras.
\end{proposition}

\subsection{Application to moduli spaces}
 Let us apply these results to  the modular operad of
 moduli spaces $\overline{\M}$.

From  th. \ref{Coperadmodformal}   it follows
\begin{corollary} $C_*(\overline{\mathcal M} ; \mathbb{Q})$ is a formal modular operad.
\end{corollary}
So,   we obtain

\begin{corollary}
Every structure of $H_*(\overline{\mathcal M} ; \mathbb{Q})$-algebra lifts to a structure
of sh $C_*(\overline{\mathcal M} ; \mathbb{Q})$-algebra.
\end{corollary}

In conclusion, we see that the minimal model $H_*(\overline{\mathcal M} ;
\mathbb{Q})_{\infty}$ of the modular operad $H_*(\overline{\mathcal M} ; \mathbb{Q})$
plays an important role in the description of the sh $H_*(\overline{\mathcal M} ;
\mathbb{Q})$-algebras and therefore of the sh $C_*(\overline{\mathcal M};\mathbb
Q)$-algebras. The explicit construction of the modular operad $H_*(\overline{\mathcal M}
; \mathbb{Q})_{\infty}$, as in the proof of th. \ref{recetapaella}, would require  the
knowledge of the homology of the moduli spaces and all its relations. We think that a
motivic formulation of this minimal modular operad (see \cite{BM}) and the determination
of the basic pieces for this building (for instance, a minimal tensor generating  family
of simple motives for the smaller abelian tensor subcategory where this operad lives)
would be a nice variant of Grothendieck's ``Lego-Teichm\"{u}ller game".

\end{large}
\end{document}